\documentclass[sn-mathphys-num]{sn-jnl}
\usepackage[ps]{skak}
\usepackage{graphicx}%
\usepackage{multirow}%
\usepackage{amsmath,amssymb,amsfonts}%
\usepackage{amsthm}%
\usepackage{mathrsfs}%
\usepackage[title]{appendix}%
\usepackage{xcolor}%
\usepackage{textcomp}%
\usepackage{manyfoot}%
\usepackage{booktabs}%
\usepackage{algorithm}%
\usepackage{algorithmicx}%
\usepackage{algpseudocode}%
\usepackage{listings}%
\usepackage{dsfont}
\usepackage{amssymb}   % symboles maths (mathbb, etc.)
\usepackage{amsmath}   % align environment
\usepackage{amsthm}    % theorems environment
\usepackage{stmaryrd}  % pour les doubles crochets-\rrbracket-\rbr
\usepackage{dsfont}    % pour l'indicatrice         \mathds{1}
\usepackage{enumitem}

\usepackage{tikz}      % faire de jolis dessins

\usepackage{pgfplots}
\pgfplotsset{width=6.5cm,compat=1.9}

\usepackage{hyperref} 
\hypersetup{
    colorlinks=true,
    linkcolor=purple,
    filecolor=magenta,      
    urlcolor=cyan,
    citecolor=teal
    }

\newcommand{\N}{\mathbb{N}}
\newcommand{\Sf}{\mathbb{S}}
\newcommand{\E}{\mathbb{E}} 
\newcommand{\T}{\mathbb{T}}

\newcommand{\R}{\mathbb{R}} 

\newcommand{\Lp}{\mathbb{L}}

\newcommand{\lbr}{\llbracket}
\newcommand{\rbr}{\rrbracket}

\newcommand{\ind}{\mathds{1}}

\newcommand{\Vect}{\mathrm{Vect}}

\newcommand{\dd}{\mathrm{d}}

\newcommand{\lnorm}{\left\|}
\newcommand{\rnorm}{\right\|}

\newcommand{\CQFD}{\hfill $\square$}

\newcommand{\und}{\mbox{ and }}

\theoremstyle{plain}
\newtheorem{lemm}{Lemma}[subsection]
\newtheorem{theo}{Theorem}[section]
\newtheorem{prop}{Proposition}[subsection]
\newtheorem{coro}{Corolary}[subsection]
\newtheorem*{propf}{Proof of the proposition}

\theoremstyle{remark}
\newtheorem*{prof}{Proof}

\newcommand{\supess}{\mathop{\smash{\mathrm{esssup}}}}

\everymath{\displaystyle} 

\usepackage{anyfontsize}

\usepackage{pgfplots}
\pgfplotsset{width=6.5cm,compat=1.9}
\usepackage{hyperref}   % hyperliens /!\ en dernier /!\
\hypersetup{
    colorlinks=true,
    linkcolor=purple,
    filecolor=magenta,      
    urlcolor=cyan,
    citecolor = blue
    }

\newcommand{\mZ}{\mathcal{Z}}
\newcommand{\mD}{\mathcal{D}}
\newcommand{\mX}{\mathcal{X}}
\newcommand{\mF}{\mathcal{F}}
\newcommand{\mN}{\mathcal{N}}
\newcommand{\mC}{\mathcal{C}}
\newcommand{\mB}{\mathcal{B}}
\newcommand{\mG}{\mathcal{G}}
\newcommand{\mP}{\mathcal{P}}
\newcommand{\mT}{\mathcal{T}}

\newcommand{\mM}{\mathcal{M}}
\newcommand{\mH}{\mathcal{H}}
\newcommand{\mE}{\mathcal{E}}

\newcommand{\bV}{\mathbf{V}}

\newcommand{\ux}{\underline{x}}
\newcommand{\us}{\underline{s}}

\newcommand{\uu}{\underline{u}}
\newcommand{\uh}{\underline{h}}

\newcommand{\uy}{\underline{y}}
\newcommand{\uv}{\underline{v}}
\newcommand{\uz}{\underline{z}}
\newcommand{\uze}{\underline{\zeta}}

\newcommand{\ut}{\underline{t}}
\newcommand{\uom}{\underline{\omega}}

\newcommand{\um}{\underline{m}}
\newcommand{\uj}{\underline{j}}
\newcommand{\uk}{{\underline{k}}}
\newcommand{\ul}{\underline{\ell}}

\newcommand{\uchi}{\underline{\smash{\chi}}}

\newcommand{\col}{\mathrm{c}}

\numberwithin{equation}{section}

\newcommand{\triple}[1]{{\left\vert\kern-0.25ex\left\vert\kern-0.25ex\left\vert #1 
    \right\vert\kern-0.25ex\right\vert\kern-0.25ex\right\vert}}

\everymath{\displaystyle} 

\begin{document}

\title[Article Title]{About a nonideal Rayleigh gas mixture model}

\subtitle{Convergence of the correlation functions}


\author[1]{\fnm{Florent} \sur{Foug\`{e}res}}\email{florent.fougeres@ens.fr}

\date{2026}

\affil[1]{\orgdiv{DMA}, \orgname{\'{E}cole normale sup\'{e}rieure}, \orgaddress{\street{45 rue d'Ulm}, \city{Paris}, \postcode{75005},  \country{France}}}

\abstract{This paper introduces a grand canonical mixture model to generalize the nonideal Rayleigh gas~\cite{2016brownian} to an asymptotically infinite amount of perturbed tagged particles. This model relies precisely on grand canonical tags, to preserve symmetry in the system, contrary to~\cite{2022mixture}. We hence define and study the convergence of the correlation functions of this system in large times, linking it to the expectancy of the empirical measure of tagged and non-tagged particles, to eventually prove a law of large numbers for this dynamics. 

We extend the quantitative study to all the correlation functions, and not only the first one, exhibiting the resultant additional factors, and we also generalize the perturbation to the whole phase space, instead of considering a space-only initial perturbation.

Eventually, we fit our adaptive time cutting~\cite{fou24} to the mixture system, even improving it to get better convergence rates. }


\ackno{I would like to thank warmly my PhD directors Isabelle Gallagher and Sergio Simonella, for their tireless help and support. Thanks also to all the persons who partake in making the department a living, stimulating and cosy place.}



\maketitle  %font blabla...

\newpage
\section{Introduction: model and results} \label{sec:intro}

\subsection{Mathematical context}

The Boltzmann equation provides a statistical model for rarefied gas dynamics. This kinetic representation may be thought about as an intermediate mesoscopic step between microscopic and macroscopic scales, and might thus be used to derive fluid mechanics equations from Newton classical equations~(see for example~\cite{1991BGL,2001ns,2008euler}). This model both allows numerous quantitative simulations~\cite{2022lbm, 2023lbm} and qualitative understanding of the intrinsic statistical behaviour of fluid matter. Indeed, introducing the \emph{entropy} of the system as a Lyapunov functional of the system, Boltzmann showed that the density of particles irreversibly converges towards an equilibrium on large time scales. This equilibrium is called the Maxwellian state: the gas occupies uniformly the available space, and the velocities of the particles get distributed according to the Maxwell Gaussian distribution~(\ref{def:maxwellian}), only depending on the temperature of the system.

The rigorous derivation of this equation from microscopic Newton equations has eventually been proved mathematically in 1975 by Oscar Erasmus Lanford III~\cite{1975lanford, 1976lanford}, in the case of hard sphere interactions. The method that Lanford used is very rigid in the tracking of particles' trajectories, hindering to prove his result for time scales larger than very small times, when only about a fifth of particles have collided. The major obstruction to large time scales is the correlation occuring between colliding particles: the quantitative estimates on the system's chaoticity deteriorate over time, making it impossible to deal with recollisions of particles.

Nonetheless, in a setting close to thermodynamic equilibrium (the Maxwellian state), the statistical stability of the dynamics guarantees a certain amount of chaos over large times, allowing a very strong control of the correlations. In the Rayleigh~gas model, describing the behaviour of a small fraction of tagged particles near equilibrium~\cite{1988spohnkineq}, a proof relying on the same ideas as Lanford's yielded the derivation \emph{on large time scales} of a linear version of the Boltzmann equation, called the linear Rayleigh--Boltzmann equation. This is the work of Henk van Beijeren, Lanford, Joel Louis Lebowitz and Herbert Spohn in 1980~\cite{1980linear}, later completed with applications to color-changing boundary conditions~\cite{1982lebspohn}.

A few decades later, in 2013, Isabelle Gallagher, Laure Saint-Raymond and Benjamin Texier reopened the work of  Lanford and his former student, Francis Gordon King~\cite{1975king}, generalizing their proof to the case of compactly supported potentials, and providing in parallel precise estimates on the convergence rate in Lanford's theorem. Eventually, these estimates led in 2016 to an article by Thierry Bodineau, Gallagher and Saint-Raymond on the convergence rate in the linear case, and its dependence on the long time scaling, so as to infer Brownian hydrodynamic limits~\cite{2016brownian}. These estimates on the convergence rates have been improved in~2024 in our previous paper~\cite{fou24}, using the adaptive time cutting that is also used in the present work.

In the present work, we introduce a generalized model for the nonideal Rayleigh~gas, where the number of tagged particles is not fixed anymore and goes to infinity with the total number of particles, yet in a small fraction so as to provide the linear Rayleigh--Boltzmann equation at the limit. This model relies on a grand-canonical mixture model, and is necessary to define the \emph{empirical measure of the tagged particle}, whose first order convergence is given in Corollary~\ref{cor:LLN}.
We study the first properties of our new Rayleigh~gas mixture model on large time scales, providing quantitative convergence results for its correlation functions in Theorem~\ref{theo1}. We extend the existing results to all the correlation functions, and not only the first one, exhibiting the big combinatorial factors that appear for high correlation functions, stemming from the time cutting method that we use. To perform this study, we harness the adaptive time cutting introduced in~\cite{fou24} to improve the quantitative convergence rates, even improving it to gain on the power of the time of validity. Performing an analysis on large time scales imposes a condition on the tagged particles' number, which is merely a technical threshold and is not required on small times.

Further statistical refinements of convergence are studied in the companion paper~\cite{fou26cumulants} thanks to the study of its~\emph{cumulants}, exploring in particular the phase transition around the low-density scaling. 

This first introducing Section~\ref{sec:intro} defines our Rayleigh~gas mixture model, with the grand canonical ensemble for tagged particles, and the main results are stated in Section~\ref{sec:main_results}. The proof of the convergence in large times of the correlation functions is given in Section~\ref{sec:proof}, using the adaptive time cutting method of~\cite{2016brownian, fou24}.

Two appendixes follow to expose technical results, about the cumulants of the exclusion (Appendix~\ref{app:cumulants_exclusion}), and the partition function combinatorics (Appendix~\ref{app:partition_function}).

\subsection{Microscopic hard sphere model} \label{sec:model2}

Microscopically, we consider exactly the hard sphere model, which resembles a perfectly elastic $d$-dimensional billiards. The state of a gas of $N$ particles is completely determined by the positions (in the $d$-dimensional torus $\T^d$) and the velocities of every particle, represented by the vector
\begin{equation*}
\uz_N = (z_1, \dots, z_N) \doteq (\ux_N, \uv_N) \in \mD^N \doteq (\T^{d} \times \R^{d})^N.
\end{equation*}
The \emph{hard sphere} model consists in an exclusion condition, which states that two particles cannot get closer than a certain diameter~$\varepsilon > 0$: the positions have to belong to the \emph{hard sphere exclusion set}
\begin{equation} \label{def:dom}
\mathcal{X}^\varepsilon_N \doteq \{ \ux_N \in \T^{dN};\ \forall\ i\neq j,\ d(x_i,x_j) > \varepsilon \},
\end{equation}
where $d(\cdot, \cdot)$ denotes the distance on the torus; and hence the state of the gas~$\uz_N$ must belong to the open domain~$\mD^\varepsilon_N \doteq \mX^\varepsilon_N \times \R^{dN}$.

 Within this set, the particles' dynamics is given by the Newton equations for uniform line movement, while on the boundary of~$\mathcal{D}^\varepsilon_N$, for at least two particles (let us say $i$ and $j$) we have that $|x_i - x_j| = \varepsilon$: they collide. Then, if the scalar product~$(x_i - x_j) \cdot (v_i - v_j)$ is positive, it means that the particles are exiting the collision in uniform line movement, but otherwise they are entering the collision and must scatter according to the following system giving the post-collisional velocities~$({v_i}', {v_j}')$:
\begin{equation} \label{eq:coll}
\left\{ \begin{array}{l}
{v_i}' = v_i - \left\langle v_i - v_j, \frac{x_i - x_j}{\varepsilon}\right \rangle\frac{x_i - x_j}{\varepsilon} \\
\ \\
{v_j}' = v_j + \left\langle v_i - v_j, \frac{x_i - x_j}{\varepsilon}\right \rangle\frac{x_i - x_j}{\varepsilon}\cdotp
\end{array} \right.
\end{equation}
In the hard sphere case, the interaction is taken instantaneous and elastic: the system~(\ref{eq:coll}) stems from the preservation of momentum and kinetic energy. 
We choose to parametrize this scattering system with the following angle
\begin{equation} \label{eq:scattering_angles}
\omega = \frac{x_i - x_j}{\varepsilon} \in \Sf^{d-1} \cdotp
\end{equation} 
This dynamics is well-defined up to a zero measure set of initial configurations, in which infinite amounts of collisions might happen in finite times, along with collisions between more than two particles at a time. This result was proved by~Roger Keith Alexander~\cite{1975alexander}, and might also be found in~\cite{2013newton}.
Some other models use non-instantaneous scattering governed by potentials of interaction that can be short-range~\cite{2013newton} or long-range~\cite{1999BoltLongrange, 2017linnoncutoff}. The review by C\'{e}dric Villani~\cite{2002villanoche} gives a global overview of collisional kinetic theory.

\subsection{Statistical description and linear Rayleigh--Boltzmann equation}
\label{sec:stat_descr_RB}
As the number of particles~$N$ gets large, this hard sphere dynamics becomes very difficult to compute, especially because it is very chaotic, so that we choose to describe the gas statistically. This paper is dedicated to the study of a gas of identical particles, yet divided in two distinguishable parts, represented by tags~$\ul_N = (\ell_1, \dots, \ell_N) \in \{0,1\}^N$: the tag~$\ell = 0$ will be attributed to particles initially distributed at thermodynamic equilibrium, and the tag~$\ell = 1$ to `tagged' particles initially perturbed from that equilibrium.

At fixed $N \in \N$ and $\varepsilon > 0$, we consider $W_N^\varepsilon(t,\uz_N, \ul_N)$ the \emph{canonical} probability density of presence of particles with tags~$\ul_N$ on the phase space $\mD^\varepsilon_N$ at time $t\geq 0$: by exchangeability, it is invariant by permutation among particles with identical tags. The microscopic dynamics provides the Liouville transport equation for $W^\varepsilon_N$ within $\mD^\varepsilon_N$ 
\begin{equation} \label{eq:liouv}
\partial_t W^\varepsilon_N + \uv_N\cdot \nabla_{\ux_N}  W^\varepsilon_N = 0.
\end{equation}
The solutions to this equation are provided by the method of characteristics and expressed in terms of the initial distribution~$W^\varepsilon_N(0)$ and of the free transport going \emph{back in time}, in such a way that we need the following boundary condition to pursue transporting back the density when two particles emerge from a collision:
\begin{equation}
d(x_i,x_j)=\varepsilon \und \langle x_i - x_j , v_i - v_j \rangle > 0 \ \ \ \ \Rightarrow\ \ \ \  W_N^\varepsilon(\uz_N) \doteq W_N^\varepsilon(\uz_N^\star), \label{eq:boundcond}
\end{equation}
where $\uz_N^\star = (z_1, \dots, x_i, {v_i}^\star, \dots,  x_j, {v_j}^\star, \dots, z_N)$ denotes the \emph{pre}-collisional state associated to $\uz_N$. Note that by reversibility, the formulas for $(v_i^\star, v_j^\star)$ are the same as for~$(v_i', v_j')$.

In the next section, we will present the grand canonical model, which randomizes the number~$N$ of particles to a random variable~$\mN$, the expectancy of which is tuned by a parameter~$\mu$ called the \emph{chemical potential}.
The kinetic limit that we consider is called the \emph{low density limit}, or Boltzmann--Grad limit, and consists in letting this chemical potential~$\mu$ go to infinity while keeping a constant mean free path $\mu^{-1} \varepsilon^{1-d} = 1$, so that the particles' diameter~$\varepsilon$ goes to 0.
In this limit, assuming initial chaos, the first marginal of the density usually converges to the solution to the Boltzmann equation~\cite{1975lanford, 2013newton, 2025DHM}.

This solution to the Boltzmann equation, studied by Boltzmann and Maxwell, when well defined, relaxes in large times to an equilibrium called the \emph{Maxwell state} and defined~\cite{1872boltzmann} as
\begin{equation} \label{def:maxwellian}
M_\beta(x,v) \doteq \left(\frac{\beta}{2\pi} \right)^{d/2} \exp\left( -\frac{\beta}{2} |v|^2 \right).
\end{equation}
The parameter~$\beta$ stands for the inverse temperature of the system, tuning its intensive (kinetic) energy. It appears that the density $M_\beta^{\otimes N} \ind_{\mD^\varepsilon_N}$ is an equilibrium of the microscopic hard sphere dynamics, and in this paper we consider specific initial conditions that are close to this thermodynamic equilibrium, to retrieve a \emph{linear version} of the Boltzmann equation, whose theory is much simpler and hence might be derived for long time scales. More precisely, we consider the \emph{nonideal Rayleigh gas} model~\cite{1988spohnkineq}: a subset of particles are tagged, breaking the particles' exchangeability. In the case of a single tagged particle, the derivation of this linear equation has been shown~\cite{2016brownian} when choosing the following perturbation of equilibrium as initial state
\begin{equation*}
W^\varepsilon_N(0, \uz_N) = \frac{\ind_{\mX^\varepsilon_N}(\ux_N)}{\mZ_N^{\varepsilon,\col}} \rho(x_1) M_\beta^{\otimes N}(\uz_N),
\end{equation*}
for $\rho \in \mC(\T^d)$ a continuous space perturbation on the torus, and $\mZ_N^{\varepsilon,\col}$ a normalization constant. Indeed, in the Boltzmann--Grad limit, the first marginal of $W^\varepsilon_N$ behaves like the solution $g\doteq M_\beta \varphi$ to the linear \emph{Rayleigh--Boltzmann equation}~\cite{2016brownian} with initial condition~$\rho$:
\begin{equation} \label{eq:phi} \def\arraystretch{1.7}
\left\{\begin{array}{rcl}
\partial_t \varphi + v\cdot \nabla_x \varphi & =&  \int_{\Sf^{d-1}} \int_{\R^d} [\varphi(v^\star) - \varphi(v)] M_\beta(v_c) \langle \omega, v_c - v \rangle_+ \dd v_c \dd \omega, \\
\varphi(0,x,v) &=& \rho(x). \end{array} \right.
\end{equation}
This linear equation~(\ref{eq:phi}) is globally well-posed in the velocity-weighted space $\Lp^\infty_x \Lp^\infty_v(M_{\beta/2})$ (see the companion paper~\cite{fou26cumulants}), and allows to derive the linear heat equation in the hydrodynamic limit~\cite{2016brownian}.

Some partial results exist for this same model with long-range interactions instead of hard sphere collisions~\cite{1999BoltLongrange, 2017linnoncutoff}, yet the complete derivation of the Rayleigh--Boltzmann~equation for general potentials is still an open problem. Other ways to derive the linear Rayleigh--Boltzmann~equation for long time scales are the \emph{ideal} Rayleigh gas model, in which the particles at equilibrium do not interact among themselves~\cite{1992centralray, 2018rayleigh, 2019rayleighannil}, and the Lorentz gas model, which consists in letting a tagged particle evolve in a frozen random background~\cite{1991spohnintpart, 1994CIP, 2008golselorentz}.

\subsection{Grand canonical framework for the mixture} \label{sec:grand_can}

As introduced in the previous section, we want to study the behaviour of an arbitrarily large subset of tagged particles perturbed away from equilibrium, whose density will follow a linear version of the Boltzmann equation. This derivation has been studied in the case of a finite set~\cite{2016brownian, fou24}, but to provide extended statistical results on this gas we hereafter take its size to infinity, yet remaining a tiny fraction of the gas. To \emph{preserve a symmetric structure} on the objets that we consider, we work in the grand canonical ensemble and introduce an additional tagging variable indicating to which set each particle belongs. This approach to describe a gas mixture is different from the canonical one used by Ioakeim Ampatzoglou, Joseph K. Miller and Nata\v{s}a Pavlovi\'{c} in their article deriving a mixed Boltzmann equation~\cite{2022mixture}; indeed their description is made in the canonical ensemble, at fixed numbers of particles of each kind.

The particles at equilibrium are taken in the usual low density (Boltzamnn--Grad) limit, whereas the tagged perturbed particles only occupy a tiny fraction of the gas, smaller than the Boltzmann--Grad density: otherwise indeed they would behave like a classical Boltzmann dilute gas, satisfying the non-linear Boltzmann equation. This all boils down to the following scaling,
\begin{equation}
  \tag{$S_{\varepsilon, \mu, \lambda}$}
  \mu \varepsilon^{d-1} = 1 \mbox{\ \ and \ \ } 1 \ll \lambda \ll \mu,
\label{eq:scaling}
\end{equation}
where~$\mu > 0$ corresponds to the chemical potential of the particles at equilibrium, and $\lambda > 0$ to that of tagged particles. The notation~$\lambda \ll \mu$ simply means that $\lambda \mu^{-1}$ goes to 0 as $\lambda$ and $\mu$ go to infinity. Formally, \textbf{the particles at initial equilibrium will be tagged with a~0, and the initially perturbed `tagged' particles will be tagged with a~1}. The tags of all particles hence form a vector $\ul_n \in \Lambda_n \doteq \{0,1\}^n$, identified to the corresponding subset $ \ul_n \subset \lbr 1, n \rbr$, with the following notation 
\begin{equation*}
|\ul_{n}| = \| \ul_n \|_1 = \left\lvert \{ i \leq n \ , \ \ell_i = 1 \} \right\rvert  \mbox{\ \ and \ \ } \varphi_0^{\otimes \ul_{n}} (\uz_{\ul_{n}}) = \prod_{\substack{i \leqslant n \\ \ell_i = 1}} \varphi_0(z_i).
\end{equation*}
Moreover, instead of $\rho$ a perturbation happening in space only, as in~\cite{2016brownian}, we will hereafter consider an initial perturbation~$\varphi_0 \in \Lp^\infty_x \Lp^\infty_v(M_{\beta/2})$ also happening according to the velocities.
Now, the mixed grand canonical ensemble consists in relaxing the number of particles, weighting it with a mixed Poisson law depending on the number of tagged particles. More precisely, at fixed~$\lambda$, and $\mu = \varepsilon^{1-d}$, we take the following weighted canonical initial densities 
\begin{equation} \label{def:initial_Wn}
W^{\varepsilon}_n(0, \uz_n, \ul_n) \doteq \frac{\lambda^{|\ul_{n}|} \mu^{n - |\ul_{n}|}}{\mZ_{\mu}} M_\beta^{\otimes n}(\uv_n) \varphi_0^{\otimes \ul_n}(\uz_{\ul_n}) \ind_{\mX_n^\varepsilon}(\ux_n)
\end{equation}
driven by the Liouville equation~\eqref{eq:liouv}, for a normalizing constant~$\mZ_\mu$ that we will adjust soon~\eqref{def:partition_function_Z}. We introduce the following \emph{correlation functions}, based on the marginals of the canonical densities. Note that we also take the marginal according to the tags, which corresponds to a sum over~$\ul^*_p \in \Lambda_p$, and we denote $\tilde{\ul}_{n+p} \doteq (\ul_n, \ul^*_p)$,  $\tilde{\uz}_{n+p} \doteq (\uz_n, \uz^*_p)$:
\begin{align*} 
F_n^\varepsilon(t, \uz_n, \ul_n)&  = \frac{1}{\mu^{n - |\ul_n|} \lambda^{|\ul_n|} } \sum_{p\geqslant 0} \sum_{\ul^*_p \in \Lambda_p } \frac{1}{p!} \int W_{n+p}^\varepsilon(t,\tilde{\uz}_{n+p}, \tilde{\ul}_{n+p})  \dd \uz^*_{p}\\
& \doteq \frac{1}{\mu^{n - |\ul_n|} \lambda^{|\ul_n|} } \sum_{p\geqslant 0} \frac{1}{p!} W_{n+p}^{\varepsilon, (n)}(t, \uz_n, \ul_n),
\end{align*}
where the normalizing (grand canonical) partition function is hence defined as follows 
\begin{align} \label{def:partition_function_Z}
\mZ_{\mu} & \doteq \sum_{p\geqslant 0} \sum_{\ul_p \in \Lambda_p } \frac{\lambda^{|\ul_{p}|} \mu^{p - |\ul_{p}|} }{p!} \int M_\beta^{\otimes p}(\uv_{p}) \varphi_0^{\otimes \ul_{p}}(\uz_{\ul_{p}})  \ind_{\mX_{p}^\varepsilon}(\ux_{p}) \dd \uz_{p}.
\end{align}
In the previous paper~\cite{2016brownian} presenting the situation where only one particle was perturbed, thanks to the invariance by translation of the system, the partition function was not depending on the perturbation$~\varphi_0$. Here, the perturbed particles are correlated one with another, preventing us from using the same argument, which is why the initial perturbation appears in the formula above. Note that the velocities of the non-perturbed particles (in the complementary of~$\ul_p$) may be integrated using the fact that the equilibrium~$M_\beta$ is of integral~1. We still keep this formulation for symmetry reasons; a more precise study of this object is made in Appendix~\ref{app:partition_function}.

Our probabilistic study is based on the random variables $(Z_{\varepsilon, i}^{[t]}, L_i)_{1\leqslant i \leqslant \mN}$, giving the states at time~$t$, and the tags of the particles. The probability density of the initial state $(Z_{\varepsilon, i}^{[0]}, L_i)$ is given above, and the evolution at time~$t$ is a deterministic piecewise affine function of the initial state. 
The correlation functions are in fact defined such that for any observable $H_n \in \mC^\infty_{\mathrm{c}}(\mD^n \times \Lambda_n)$, we have
\begin{align}
\E\left[ \sum_{1 \leq i_k \neq i_j \leq \mathcal{N}} H_n(Z_{\varepsilon,i_1}^{[t]}, L_{i_1}, \dots, Z_{\varepsilon, i_n}^{[t]}, L_{i_n})  \right] 
& = \E\left[ \delta_{\mathcal{N} \geq n} \frac{\mathcal{N}!}{(\mathcal{N}-n)!} H_n(Z_{\varepsilon,1}^{[t]},L_1, \dots, Z_{\varepsilon,n}^{[t]},L_n)\right] \nonumber \\
\ & = \sum_{p = n}^\infty \frac{1}{p!} \frac{p!}{(p-n)!} \sum_{\ul_p \in \Lambda_p} \int_{\mathcal{D}^\varepsilon_p} W^\varepsilon_{p}(t,\uz_p, \ul_p) H_n(\uz_n, \ul_n) \dd \uz_p \nonumber \\
\ & = \sum_{\ul_n  \in \Lambda_n} \mu^{n - |\ul_n|} \lambda^{|\ul_n|}  \int_{\mathcal{D}^\varepsilon_n} F_n^\varepsilon(t,\uz_n, \ul_n)  H_n(\uz_n, \ul_n) \dd \uz_n. \label{eq:empmes}
\end{align} 
By the usual BBGKY arguments~\cite{fouthese, 2013newton} the canonical densities~(\ref{def:initial_Wn}) satisfy the following hierarchy
\begin{align}
\partial_t W_N^{\varepsilon, (n)} + \uv_n \cdot \nabla_{\ux_n} W_N^{\varepsilon, (n)} & = (N-n)\varepsilon^{d-1} \sum_{\ell_{}= 0}^1 \sum_{i=1}^n \int \langle \omega, v_{\col} - v_i \rangle W_N^{\varepsilon, (n+1)}( \uz_n, \ul_n, x_i + \varepsilon\omega, v_{\col}, \ell_{})\dd \omega \dd v_{\col}\nonumber \\
& \doteq (N-n)\varepsilon^{d-1} \left( \mC^{\langle 0 \rangle}_n  W^{\varepsilon, (n+1)}_N + \mC^{\langle 1 \rangle}_n W_N^{\varepsilon, (n+1)} \right). \label{def:collision_operator}
\end{align}
The two operators we hence define are called \emph{collision operators} and represent the action of a $(n+1)$-th particle (tagged~0 or tagged~1) on the average distribution of $n$~particles. To better understand these operators, we rewrite them using the boundary condition~\eqref{eq:boundcond} when $\langle \omega, v_{n+1} - v_i \rangle > 0$, and the change of variable $\omega \mapsto - \omega$ otherwise:
\begin{align} 
\mC_n^\ell W_N^{\varepsilon, (n+1)} = \sum_{i=1}^n  \int &\dd \omega \dd v_{n+1} \langle \omega, v_{n+1} - v_i \rangle_+ \times \nonumber \\ 
&  \Bigl[W_N^{\varepsilon, (n+1)}(\uz_n^\star, \ul_n, x_i + \varepsilon \omega, v_{n+1}^\star, \ell) - W_N^{\varepsilon, (n+1)}(\uz_n, \ul_n, x_i - \varepsilon \omega, v_{n+1}, \ell) \Bigr] . \label{def:collision_op2}
\end{align}
Morally, we look at the influence of a $(n+1)$-th particle---with tag~$\ell$---on the dynamics, colliding with one of the $n$ existing ones with angle~$\omega$ and velocity~$v_{n+1}$, whence the name collision operators.
 The \emph{cross section} $\langle \omega, v_{n+1} - v_i\rangle_+$ weights the likelihood of such a collision. 
As a consequence, the correlation functions satisfy 
\begin{align*}
\partial_t F_n^\varepsilon + \uv_n \cdot \nabla_{\ux_n} F_n^\varepsilon & = \frac{1}{\mu^{n-{|\ul_n|}} \lambda^{|\ul_n|}} \sum_{k\geqslant 0} \frac{k \varepsilon^{d-1}}{k!} \left( \mC^{\langle 0 \rangle}_n  W_{n+k}^{\varepsilon, (n+1)} + \mC^{\langle 1 \rangle}_n W_{n+k}^{\varepsilon, (n+1)}\right)\\
& = \frac{\mu \varepsilon^{d-1}}{\mu^{n+1-{|\ul_n|}} \lambda^{|\ul_n|}} \sum_{k\geqslant 0} \frac{1}{k!}\mC^{\langle 0 \rangle}_n  W_{n+k+1}^{\varepsilon, (n+1)} + \frac{\lambda \varepsilon^{d-1}}{\mu^{n-{|\ul_n|}} \lambda^{|\ul_n|+1}} \sum_{k\geqslant 0} \frac{ \varepsilon^{d-1}}{k!} \mC^{\langle 1 \rangle}_n W_{n+k+1}^{\varepsilon, (n+1)},
\end{align*}
which eventually leads to the following hierarchy in our mixed Boltzmann--Grad scaling~\eqref{eq:scaling}
\begin{equation} \label{eq:BBBGKY}
\partial_t F_n^\varepsilon + \uv_n \cdot \nabla_{\ux_n} F_n^\varepsilon = \mC^{\langle 0 \rangle}_n  F^\varepsilon_{n+1} + \frac{\lambda}{\mu} \mC^{\langle 1 \rangle}_n F^\varepsilon_{n+1}.
\end{equation}
The formal limit of the hierarchy above~\eqref{eq:BBBGKY} with initial conditions~\eqref{def:initial_Wn}, presented in Section~\ref{sec:strategy}, is satisfied by the family~$(M_\beta^{\otimes n} \varphi^{\otimes \ell_n})_{n \geqslant 1}$, where~$\varphi$ is the solution of the Rayleigh--Boltzmann equation~\eqref{eq:phi} with initial data~$\varphi_0$. This result is formalized in the following section, with a quantitative convergence~rate.

\subsection{Main results: long-time convergence of the correlation functions} \label{sec:main_results}

The following theorem provides a convergence rate of the mixed correlation functions to the solutions of the Rayleigh--Boltzmann equation, which is a generalization of~\cite{fou24, 2016brownian} with a time scale of validity improved by a power $1/4$ thanks to a more precise computation (see the proof of Proposition~\ref{prop:estrem}). We extend the result to all the correlation functions (not only the first one), yet at the cost of a bad constant~$n^{cn}$ stemming from the time cutting method we use. This constant, which did not appear in the convergence of the first marginal, is due to an accumulation of errors at each time step of our cutting. Note finally that we use here the adaptive time cutting introduced in~\cite{fou24}, improving greatly the convergence rate compared with~\cite{2016brownian}.
\begin{theo}[Convergence of the correlation functions] \label{theo1}
For some sets~$\Delta_n^\varepsilon \subset \mD^n$ whose measure goes to 0 with~$\varepsilon$, there exists a constant $c_\beta$ depending only on the temperature and the dimension such that, for any $\alpha \in (0,3/4)$, as long as 
\begin{equation} \label{ineq:condt} t \lesssim \left( \log\left| c_\beta \log \varepsilon \right|\right) ^{\frac{3}{4} - \alpha} \ \ \und \ \ \lambda \lesssim |\log \varepsilon|^{1 - \alpha},
\end{equation}
and for $\varepsilon$ small enough, one has the following convergence rate of the correlation functions to the linear Rayleigh--Boltzmann solutions in our mixed low density scaling~\eqref{eq:scaling}, for a constant~$c > 0$;
\begin{equation*}
\| F_n^\varepsilon - M_\beta^{\otimes n} \varphi^{\otimes \ul_n}\|_{\Lp^{\infty}([0,t] \times \mD^n \setminus \Delta_n^\varepsilon)} \leq n^{cn} \exp\left( - c_\beta |\log \varepsilon |^{1- \alpha}  \right).
\end{equation*}
\end{theo}
The notation $t \lesssim \left( \log\left| c_\beta \log \varepsilon \right|\right) ^{\frac{3}{4} - \alpha}$ means that, for a good constant $c>0$ depending only on the dimension~$d$ and the inverse temperature~$\beta$, one has 
\begin{equation*}
t \leq c \left( \log\left| c_\beta \log \varepsilon \right|\right) ^{\frac{3}{4} - \alpha}.
\end{equation*}

The proof of this theorem is the subject of Section~\ref{sec:proof}. It is close to the proof presented in~\cite{fou24}, but in the mixed grand canonical framework and for all the correlation functions instead of the first marginal only.

This theorem provides a first corollary on the statistical behaviour of the gas. So as to state this result, for any observable~$H \in \mC^\infty_{\mathrm{c}}(\mD \times \Lambda_1)$ we define the following random variables: the \emph{empirical measure of all particles}
\begin{equation} \label{def:emp_meas}
\pi^\varepsilon_t[H] \doteq \frac{1}{\mu} \sum_{i=1}^\mN H(Z^{[t]}_{\varepsilon,i}	, L_i),
\end{equation}
and the \emph{empirical measure of tagged particles}
\begin{equation} \label{def:emp_meas_tag}
\tilde{\pi}^\varepsilon_t[H] \doteq \frac{1}{\lambda} \sum_{i=1}^\mN H(Z^{[t]}_{\varepsilon,i}, L_i) \ind_{L_i = 1}.
\end{equation}
Thanks to Theorem~\ref{theo1}, one can deduce the convergence of these empirical measures, hence providing a law of large numbers for the hard sphere dynamics. This result is given in the space~$\mathbf{L}^2$ of \emph{square-integrable random variables}, endowed with the norm~$\E\left[|\cdot|^2\right]$.
\begin{coro}[Law of large numbers for the dynamics] \label{cor:LLN}
The empirical measures converge as random variables in~$\mathbf{L}^2$, the non-tagged particles towards an equilibrium state, and the tagged ones to a state described by the linear Rayleigh--Boltzmann equation~\eqref{eq:phi}, in the following way
\begin{equation} \label{eq:LLN0}
\pi^\varepsilon_t[H] \xrightarrow[\varepsilon \to 0]{\mathbf{L}^2} \int M_\beta(v)  H(z, 0) \dd z,
\end{equation}
and 
\begin{equation} \label{eq:LLN}
\tilde{\pi}^\varepsilon_t[H] \xrightarrow[\varepsilon \to 0]{\mathbf{L}^2} \int M_\beta(v) \varphi(t,z) H(z, 1) \dd z.
\end{equation}
\end{coro} 
\begin{prof} To show that the random variable $\pi^\varepsilon_t[H]$ converges in $\mathbf{L}^2$ to a deterministic limit $a \in \R$, writing 
\begin{equation*}
\pi^\varepsilon_t[H] - a = \pi^\varepsilon_t[H] - \E\bigl[ \pi^\varepsilon_t[H] \bigr] + \E\bigl[ \pi^\varepsilon_t[H] \bigr] - a ,
\end{equation*}
it is enough to show that $\E\bigl[ \pi^\varepsilon_t[H] \bigr] \xrightarrow[\varepsilon \to 0]{} a $ and $\E\left[ \bigl\lvert \pi^\varepsilon_t[H] - \E\bigl[ \pi^\varepsilon_t[H] \bigr] \bigr\rvert^2 \right] \xrightarrow[\varepsilon \to 0]{} 0$.
Using formula~\eqref{eq:empmes}, one can write
\begin{align*}
\E\left[ \frac{1}{\mu} \sum_{i=1}^\mN H(Z_{\varepsilon,i}^{[t]}, L_i) \right] & = \int F_1^\varepsilon(t,z_1, 0) H(z_1, 0) \dd z_1 + \frac{\lambda}{\mu} \int F_1^\varepsilon(t,z_1, 1) H(z_1, 1) \dd z_1,
\end{align*} 
so that the expectancies converge by Theorem~\ref{theo1} above. For concision, we show the  fact that the variance vanishes in the tagged case, denoting $h \doteq H(\cdot, 1)$. The equilibrium case is treated in a similar though simpler way. Let us compute, once again by formula~\eqref{eq:empmes},
\begin{align*}
\E\left[ \bigl| \tilde{\pi}^\varepsilon_t - \E\bigl[ \tilde{\pi}^\varepsilon_t[H] \bigr] \bigr|^2 \right] & = \frac{1}{\lambda^2} \E\left[  \sum_{i,j=1}^n h(Z_{\varepsilon,i}^{[t]})h(Z_{\varepsilon,j}^{[t]} ) \right] -  \frac{2}{\lambda} \E\left[ \sum_{i=1}^n h(Z_{\varepsilon,i}^{[t]} ) \int F_1^\varepsilon(1) h \right] + \left(\int F_1^\varepsilon(1) h\right)^2\\
& = \int F_2^\varepsilon(1,1) h^{\otimes 2} + \frac{1}{\lambda} \int F_1^\varepsilon(1) h^2 - 2 \left(\int F_1^\varepsilon(1) h\right)^2 + \left(\int F_1^\varepsilon(1) h\right)^2 \\
& \xrightarrow[\varepsilon \to 0]{} \int M_\beta^{\otimes 2} \varphi^{\otimes 2} + 0 - \left( \int M_\beta \varphi \right)^2 = 0,
\end{align*}
thanks to the convergence of the correlation functions (Theorem~\ref{theo1}).
The other convergence~\eqref{eq:LLN0}  follows similarly in our scaling, using the fact that $\frac{\lambda}{\mu}$ goes to~0. \CQFD
\end{prof}

\section{Proof of the convergence of the correlation functions} \label{sec:proof}

We give in this section the proof of Theorem~\ref{theo1}.

\subsection{Pseudo-trajectories and strategy of proof} \label{sec:strategy}

Let us choose to denote $p_\mu \doteq \frac{\lambda}{\mu}$ the fraction of initially perturbed particles, so that iterating Duhamel formula as in~\cite{2013newton} or~\cite{2016brownian}, we can write the Dyson~expansion
\begin{equation} \label{eq:Dyson}
F_n^\varepsilon(t) = \sum_{k\geqslant 0} \sum_{\ul^*_k \in \Lambda_k } p_\mu^{|\ul^*_k|} Q_{n, \ul^*_k}(t) F^\varepsilon_{n+k}(0),
\end{equation}
developing the choice of the encountered tags~$\ul^*_k \doteq (\tilde{\ell}_{n+1}, \dots, \tilde{\ell}_{n+k})$, with the successive-collision operators defined as
\begin{equation} \label{def:successive_collision_operator}
Q_{n, \ul^*_k}(t) \doteq \int_{T_k(t)} \Theta_n(t-t_1)\mathcal{C}^{\tilde{\ell}_{n+1}}_n \Theta_{n+1}(t_1 - t_2)\dots \mC^{\tilde{\ell}_{n+k}}_{n+k-1} \Theta_{n+k}(t_k)  \dd \ut_{k},
\end{equation}
where $\Theta_n(\tau)$ denotes the transport semi-group operator in $\mD^\varepsilon_n$ with specular reflections, for a time~$\tau$. The collision times are integrated over
\begin{equation} \label{def:time_ensemble}
T_k(t) \doteq \Bigl\{\ \ut_k \ \Bigl | \ 0 \doteq t_{k+1} \leq t_k \leq \dots \leq t_1 \leq  t_0 \doteq t \Bigr\}.
\end{equation}
The main idea of the proof, coming from Lanford's original paper~\cite{1975lanford}, is to use a coupling between this expansion and its limit version, implying imaginary histories of the particles, among those that eventually lead to the state $\uz_n$ at time~$t$. These histories, called \emph{pseudo-trajectories}, are non-physical trajectories that---in a way---allow to extend the method of characteristics for the successive-collision operators. 

Indeed, the transport operators appearing in~\eqref{def:successive_collision_operator} correspond to following the characteristics of free transport, with specular reflections: taking the first operator $\Theta_n(t-t_1)$ of a functional amounts to considering this functional at time $t_1$,  in a state~$\uz_n^{[t_1]}$ given by the backwards hard sphere dynamics. 

Then, the first collision operator~\eqref{def:collision_op2} writes
\begin{align*} 
&\mC_n^\ell F_{n+1}^\varepsilon = \sum_{i=1}^n \sum_{s_{1} = \pm 1} s_{1} \int \dd \omega_1 \dd v_{n+1} \langle \omega_1, v_{n+1} - v_i \rangle_+  F_{n+1}^\varepsilon(\uz_n^{\langle s_{1}\rangle}, \ul_n, x_i + s_1 \varepsilon \omega_1, v_{n+1}^{\langle s_{1}\rangle}, \ell),
\end{align*}
where $\uz_n^{\langle +1 \rangle} = \uz_n^{\star}$ and $\uz_n^{\langle -1 \rangle} = \uz_n$, scattered for the gain term, and let unchanged for the loss term, so that the collision is always incoming, allowing to pursue the backwards method of characteristics with the next transport operator. Hence, for given collision parameters~$(i, s_1, \omega_1, v_{n+1})$, this operator can be seen as a weighted adjunction of a particle to the characteristics---or pseudo-trajectory---which scatters (or not, according to $s_1$) with particle~$i$, creating a new state~$\uz_{n+1}^{[t_1]} \doteq (\uz_{n}^{[t_1]}, x_i+s_1 \varepsilon \omega_1, v_{n+1}^{\langle s_1 \rangle})$. The integration and sum over these collision parameters will yield an integral over pseudo-trajectories. Iterating this extended method of characteristics and tracking the pseudo-trajectories $(\uz_{n+j}^{[t_j]})$ thus constructed, we bring the analysis back to the value of the functional at time~$\tau = 0$, in the state~$\uz_{n+k}^{[0]}$. 

We will have to record the numbering labels of the existing particles meeting the new one, the velocities of the particles that spring up, the angles at which the encounters happen, and whether they scatter or not.
The pseudo-trajectories will also keep track of the tags of the encountered particles.

Here is precisely how we construct the pseudo-trajectories.
The choice of the successive encountered tags is registered in $\ul^*_k = (\tilde{\ell}_{n+1}, \dots, \tilde{\ell}_{n+k})$, and expanding all the sums in all the collision operators~(\ref{def:successive_collision_operator}), we can sum them up to the history~$(m_1, \dots, m_k)$ of which particle encountered the $(n+i)$-th new one. These particles naturally belong to the following set 
\begin{equation*} \label{eq:Mnk}
\mM_{n,k} \doteq \Bigl\{ (m_1, \dots, m_k) \ \Bigl| \ \forall i \leq k,  m_i \leq n+i-1 \Bigr\}.
\end{equation*}
We consider the scattering labels~$(s_1, \dots, s_k) \in \{ \pm 1\}^k$. The fact that some encounters do not scatter, along with the fact that some particles are artificially added, is why the pseudo-trajectories are not physical trajectories.
Once the total history 
\begin{equation} \label{def:history}
\uchi_k \doteq (\um_k, \ul^*_k, \us_k) \in \mH_{n,k} \doteq \mM_{n,k} \times \Lambda_k \times \{\pm 1\}^k
\end{equation}
is fixed, for given collision parameters~$(\uom_k, v_{n+1}, \dots, v_{n+k})$, we can construct the pseudo-trajectories for every endstate~$\uz_n =\uz_n^{[t]}$, backwards in time to an initial configuration~$\uz_{n+m}^{[0]}$, following the inductive procedure below:
\begin{equation} \label{eq:pseudotraj}
\left\{ \def\arraystretch{1.5} \begin{array}{l}
\uz_n^{[t]} \doteq \uz_n  \\
\forall \ i \in \lbr 0, k \rbr,\ \forall \tau \in (t_{i+1}, t_i),\ \uz_{n+i}^{[\tau]} \mbox{ follows (backwards) the physical hard sphere dynamics}\\
\forall \ i \in \lbr 1, k \rbr,\ \uz_{n+i}^{[t_{i}]} = \left(\uz_{n+i-1}^{[t_i^+],\langle s_{i} \rangle }, x_{m_i} + s_i \varepsilon \omega_i, v_{n+1}^{\langle s_{i} \rangle}\right).
\end{array} \right.
\end{equation}
One may observe that the change of velocities in the last step is automatic by the hard sphere dynamics' boundary condition, but it will not be for the limit version of pseudo-trajectories, since the limit particles are formally pointwise.
In the end, one can write the \emph{pseudo-trajectory formulation} of the Dyson expansion
\begin{equation}
 \label{eq:pseudo_traj_form}
F_{n}^\varepsilon(t) = \sum_{k \geqslant 0} \sum_{\uchi_k} p_\mu^{|\ul^*_k|} \int_{T_k(t)} \dd \ut_k \int \dd \underline{\omega}_k \dd v_{n+1} \dots \dd v_{n+k} \prod_{i=1}^{k} s_i \langle \omega_{i},v_{n+i} - v_{m_i}^{[t_i^+]}  \rangle_+ F_{n+k}^\varepsilon\bigl(0, \uz_{n+k}^{[0]}, \tilde{\ul}_{n+k} \bigr).
\end{equation} 
A small technical detail lies in the fact that the added particles must satisfy the exclusion condition. A way to deal with it may be to impose a condition on the domain of integration of the collision angles~\cite{2023grandev}, yet here to simplify we merely change the definition of the pseudo-trajectories: if at any moment the exclusion condition is violated by the adjunction of a particle, then the trajectories are frozen in this state until time~$\tau = 0$, so that the integral formally vanishes thanks to the initial distribution~$F_{n+k}^\varepsilon(0)$ being~0 outside of $\mD^\varepsilon_{n+k}$. 

Indeed, our goal is now to prove the convergence of the correlation functions to the \textbf{limit densities}
 \begin{equation*}
 G_n(t, \uz_n, \ul_n) \doteq M_\beta^{\otimes n}(\uv_n) \varphi^{\otimes \ul_n}(t,\uz_{\ul_n}), 
 \end{equation*}
where $\varphi$ is the solution to the linear Rayleigh--Boltzmann~equation~(\ref{eq:phi}). Because of the structure of this equation, this family satisfies the following hierarchy
\begin{equation} \label{eq:limit_hierarchy}
(\partial_t + \uv \cdot \nabla_{\ux}) G_n = \sum_{i=1}^n \sum_{s_\col = \pm 1} s_\col \int \dd \omega \dd v_{n+1} \langle \omega, v_{n+1} - v_i \rangle_+ G_{n+1}(\uz_n^{\langle s_\col \rangle}, \ul_n, x_i, v_{n+1}^{\langle s_\col \rangle}, 0), 
\end{equation}
noticing that the terms vanish when the scattering occur between two~particles distributed according to the equilibrium~$M_\beta^{\otimes n}$. This equation is the formal limit of the BBGKY hierarchy~\eqref{eq:BBBGKY} in the mixed low density regime~\eqref{eq:scaling}: it makes only appear the collision operator linked to equilibrum particles, tagged~0, since the other one has a factor $\lambda / \mu$ that vanishes at the limit. It leads to the following limit version of the Dyson expansion~\eqref{eq:Dyson}
\begin{equation} \label{eq:Dyson_lim}
G_n(t) = \sum_{k\geqslant 0}  Q^{\lim}_{n, \underline{0}_k}(t) G_{n+k}(0),
\end{equation}
where the following successive-collision operators contain only collisions with particles at equilibrium, and limit free transport of pointwise particles, without scattering:
\begin{equation} \label{def:successive_collision_operator_lim}
Q^{\lim}_{n, \underline{0}_k}(t) \doteq \int_{T_k(t)} \Theta^{\lim}_n(t-t_1)\mathcal{C}^{\langle 0 \rangle, \lim}_n   \Theta^{\lim}_{n+1}(t_1 - t_2)\dots \mC^{\langle 0 \rangle, \lim}_{n+k-1} \Theta^{\lim}_{n+k}(t_k)  \dd \ut_{k}.
\end{equation}
The limit collision operators~$\mathcal{C}^{\langle 0 \rangle, \lim}_n$ are the formal limit of the collision operators~\eqref{def:collision_op2}, for~$\varepsilon = 0$.
The same computation as below for this limit hierarchy leads to a similar writing in terms of pseudo-trajectories
\begin{equation}
 \label{eq:pseudo_traj_form_lim}
G_{n}(t) = \sum_{k \geqslant 0} \sum_{\uchi_k} \ind_{\ul^*_k = \underline{0}_k} \int_{T_k(t)} \dd \ut_k \int \dd \underline{\omega}_k \dd v_{n+1} \dots \dd v_{n+k} \prod_{i=1}^{k} s_i \langle \omega_{i}, v_{n+i} - v_{m_i}^{[t_i^+]} \rangle_+ G_{n+k}\bigl(0, \uze_{n+k}^{[0]}, \tilde{\ul}_{n+k} \bigr),
\end{equation} 
where the limit pseudo-trajectories~$(\uze_{n+i}^{[t_i]})_{i\leqslant k }$ are defined as their hard sphere versions~\eqref{eq:pseudotraj} for $\varepsilon = 0$, with the noticeable difference that in the dynamics followed on each time interval~$(t_{i+1}, t_i)$, the particles are pointwise and hence follow the free flow without any scattering.

\paragraph{Strategy of proof} We will couple both pseudo-trajectory formulations~\eqref{eq:pseudo_traj_form} and~\eqref{eq:pseudo_traj_form_lim}, bringing down the difference at a certain time~$(F_n^\varepsilon - G_n)(t, \uz_n)$ to the difference at time~0 of higher correlation functions $\left[F_{n+k}^\varepsilon(0, \uz_{n+k}^{[0]}) - G_{n+k}(0, \uze_{n+k}^{[0]})\right]$. For these hierarchies to be well coupled, the two pseudo-trajectories must be close; the classical argument is to show that the transport operators $\Theta_k$, appearing between the collision operators, do not imply additional recollisions, since the trajectories would diverge from the limit transport operator~$\Theta_k^{\lim}$, defined on the whole space $\mD^k$ without recollisions.

Indeed with this method, it will be enough to use \textbf{continuity estimates} on the operators to bring the convergence back to time~$\tau=0$ where we can use explicit \textbf{initial proximity}. Nevertheless, these continuity estimates demand to work with trajectories that do not contain too many particles, so that we will first compute a \textbf{tree pruning}, and control the pruned-out term using some \textbf{a priori estimates} on the densities.
The last step before proving the \textbf{convergence} will be to \textbf{discard the pseudo-trajectories} in which some perturbed particles are encountered, as this situation does not happen in the limit pseudo-trajectories. 
This strategy will be followed in the following sections.

\subsection{Initial proximity} 

The total kinetic energy, preserved at fixed number of particles by the transport and by elastic collisions, will be denoted
\begin{equation*}
\| \uv_k\|^2 \doteq \sum_{i = 1}^k |v_i|^2,
\end{equation*} 
where $|v_i|$ is the Euclidean norm of the velocity $v_i \in \R^d$ of particle~$i$. 
For an inverse temperature~$\beta > 0$ and $k\in \N^*$, we consider the space $\mF_{k, \beta}$ of measurable functions defined almost everywhere on the domain $\mD^k$ such that
\begin{equation} \label{def:foncspace}
\|f_k\|_{k, \beta} \doteq \supess_{\uz_k \in \mD^k} \bigl\lvert f_k(\uz_k) \exp( \beta \|\uv_k\|^2) \bigr \rvert < \infty,
\end{equation}
hence decreasing at least as the Gaussian equilibrium~$M_{2\beta}^{\otimes k}$ in velocities.  We denote
\begin{equation} \label{def:C0}
C_0 \doteq \max \left[ \| M_\beta \|_{1, \beta/2} \ ; \ \| M_\beta \varphi_0 \|_{1, \beta/4}\ ; \  \left\| M_\beta M_{\beta/2}^{-\frac{1}{2}} \varphi_0 \right\|_{\Lp^\infty(\mD)}  \right].
\end{equation}
The initial error between the microscopic densities and the limit ones is mainly due to the exclusion condition, of which we can compute an explicit control. 
\begin{prop}[Initial proximity] \label{prop:initial_proximity} For all~$n\in \N$ and tags~$\ul_n \in \Lambda_n$, for any $\lambda, \mu > 0$ in the scaling~\eqref{eq:scaling}, one has
\begin{equation} 
\left\| \ind_{\mX_n^\varepsilon}M_\beta^{\otimes n} \varphi_0^{\otimes \ul_n} - F_n^\varepsilon(0, \ul_n) \right\|_{n, \beta/4}\leq  C_0^n \varepsilon.
\end{equation}
\end{prop} 
\begin{propf} We denote $\ind_{i \not\sim j} \doteq \ind_{d(x_i, x_j) > \varepsilon}$ and $ \ind_{i \sim j} \doteq \ind_{d(x_i,x_j) \leq \varepsilon}$ the indicator of exclusion between~$i$ and~$j$ and its complementary. Recalling that $\ind_{\mX_N^\varepsilon}$ denotes the exclusion condition~\eqref{def:dom} and by definition~\eqref{def:partition_function_Z} of the partition function $\mZ_\mu$, we can write, once again denoting~$\tilde{\ul}_{n+p} \doteq (\ul_n, \ul^*_p)$:
\begin{align*}
& \ind_{\mX^\varepsilon_n}M_\beta^{\otimes n} \varphi_0^{\otimes \ul_n} - F_n^\varepsilon(0) \\
& \hspace{10mm} = \frac{1 }{\mZ_{\mu}} \sum_{p\geqslant 0} \frac{1}{p!} \left[ \sum_{\ul^*_p  \in \Lambda_p} \lambda^{|\ul^*_p |} \mu^{p - |\ul^*_p |} \ind_{\mX^\varepsilon_n} M_\beta^{\otimes n} \varphi_0^{\otimes \ul_n} \int \varphi_0^{\otimes \ul^*_p } M_\beta^{\otimes p} \ind_{\mX^\varepsilon_p} - \left(\varphi_0^{\otimes \tilde{\ul}_{n+p}} M_\beta^{\otimes n+p}\ind_{\mX^\varepsilon_{n+p}}\right)^{(n)}\right].
\end{align*}
Expanding the marginals' formula, the following term appears, in which we decompose the exclusion indicator~$\ind_{\mX_{n+p}^\varepsilon}$ as
\begin{align*}
& \ind_{\mX^\varepsilon_n} M_\beta^{\otimes n} \varphi_0^{\otimes \ul_n} \int \varphi_0^{\otimes \ul^*_p} M_\beta^{\otimes p} \ind_{\mX^\varepsilon_p} - \int \varphi_0^{\otimes \tilde{\ul}_{n+p}} M_\beta^{\otimes n+p}\ind_{\mX^\varepsilon_{n+p}} \dd \uz^*_p\\
& \hspace{20mm} = \ind_{\mX^\varepsilon_n} M_\beta^{\otimes n} \varphi_0^{\otimes \ul_{n}} \int  (M_\beta\varphi_0)^{\otimes \ul^*_p}(\uz^*_{\ul^*_p }) \ind_{\mX^\varepsilon_{p}}(\ux^*_{p}) \left(1 - \prod_{\substack{1 \leqslant i \leqslant n \\ 1 \leqslant j \leqslant p}} \ind_{x_i \not\sim x^*_j}  \right)  \dd \ux^*_p \uv^*_{\ul^*_p },
\end{align*}
using that $M_\beta$ is of integral~1. Then, we will harness the following basic set property,
\begin{align*}  1 - \prod_{\substack{1 \leqslant i \leqslant n \\ 1 \leqslant j \leqslant p}} \ind_{x_i \not\sim x^*_j}  \leq \sum_{\substack{1 \leqslant i \leqslant n \\ 1 \leqslant j \leqslant p}} \ind_{x_i \sim x^*_j},
\end{align*}
yielding
\begin{align*}
\left\lvert \ind_{\mX^\varepsilon_n}M_\beta^{\otimes n} \varphi_0^{\otimes \ul_n} - F_n^\varepsilon(0) \right\rvert& \leq \frac{\ind_{\mX^\varepsilon_n} M_\beta^{\otimes n} \varphi_0^{\otimes \ul_n} }{\mZ_{\mu}} \sum_{p\geqslant 0} \sum_{\ul^*_p \in \Lambda_p } \frac{\lambda^{|\ul^*_{p}|} \mu^{p - |\ul^*_{p}|} }{p!}  \sum_{\substack{1 \leqslant i \leqslant n \\ 1 \leqslant j \leqslant p}} 
\int (M_\beta\varphi_0)^{\otimes \ul^*_{p} } \ind_{\mX^\varepsilon_{p}} \ind_{x_i \sim x^*_j} \dd \ux^*_p \dd \uv_{\ul_p^*}^*.
\end{align*}
To be able to integrate over~$x^*_j$, we denote~$\check{\ul}_p^{(j)} \doteq (\ell^*_1, \dots, \ell^*_{j-1}, \ell^*_{j+1}, \dots, \ell^*_p)$ the vector of all tags apart from~$j$. Using the definition~\eqref{def:C0} of~$C_0$ and summing over~$\ell_j^*$, we get
\begin{align*}
& \left\|\ind_{\mX^\varepsilon_n}M_\beta^{\otimes n} \varphi_0^{\otimes \ul_n} - F_n^\varepsilon(0)\right\|_{n, \beta} \\
& \hspace{25mm} \leq \frac{C_0^n }{\mZ_{\mu}} \sum_{p\geqslant 0} \ \sum_{\substack{1 \leqslant i \leqslant n \\ 1 \leqslant j \leqslant p}} \ \sum_{\check{\ul}_p^{(j)} \in \Lambda_{p-1} } \frac{\lambda^{|\check{\ul}_p^{(j)}|} \mu^{p - 1 - |\check{\ul}_p^{(j)}|} }{p!} (\lambda C_0 + \mu) \int (M_\beta\varphi_0)^{\otimes \check{\ul}^{(j)}_p } \ind_{\mX^\varepsilon_{p}} \ind_{x_i \sim x^*_j} \\
& \hspace{25mm} \leq \frac{C_0^n }{\mZ_{\mu}} \sum_{p\geqslant 0} np \sum_{\ul^*_{p-1} \in \Lambda_{p-1} } \frac{\lambda^{|\ul^*_{p-1}|} \mu^{p -1 - |\ul^*_{p-1}|} }{p!} 2 \mu |\Sf^{d-1}| \varepsilon^{d} \int (M_\beta\varphi_0)^{\otimes \ul^*_{p-1} } \ind_{\mX^\varepsilon_{p-1}} 
\end{align*}
using the exchangeability of identical particles, integrating over~$x_j^*$ and using the fact that~$C_0 \lambda < \mu$. We get in the mixed low density scaling~\eqref{eq:scaling}
\begin{align*}
\left\|\ind_{\mX^\varepsilon_n}M_\beta^{\otimes n} \varphi_0^{\otimes \ul_n} - F_n^\varepsilon(0)\right\|_{n, \beta} &  \leq 2  |\Sf^{d-1}| \varepsilon \frac{C_0^n }{\mZ_{\mu}} n \sum_{p\geqslant 1}  \sum_{\ul^*_{p-1} \in \Lambda_{p-1} } \frac{\lambda^{|\ul^*_{p-1}|} \mu^{p -1 - |\ul^*_{p-1}|} }{(p-1)!}  \int (M_\beta\varphi_0)^{\otimes \ul^*_{p-1} } \ind_{\mX^\varepsilon_{p-1}},
\end{align*}
which concludes the proof recognizing the partition function~$\mZ_\mu$~\eqref{def:partition_function_Z} after an index shift.
\CQFD
\end{propf}

\subsection{A priori estimates}
\label{sec:apriori}

As in the previous works~\cite{2016brownian, fou24} about the Rayleigh gas, the long-time derivation is allowed thanks to a priori estimates yielded by the rigid structure of the equilibrium. Here, some additional technical difficulties, dealt with in Appendix~\ref{app:partition_function}, appear in the proof of these estimates because of the structure of the grand canonical mixture and its partition function. First, let us observe that the initial canonical densities defined in~\eqref{def:initial_Wn} satisfy for all $(\uz_n, \ul_n) \in \mD^n \times \Lambda_n$ that
\begin{equation} \label{ineq:W_n0}
W_n^\varepsilon(0, \uz_n, \ul_n) \leq \frac{\lambda^{|\ul_{n}|} \mu^{n - |\ul_{n}|}}{\mZ_{\mu}} \left\| M_\beta M_{\beta/2}^{-\frac{1}{2}} \varphi_0 \right\|_{\Lp^\infty}^{|\ul_n|} M_{\beta/2}^{\otimes n}(\uv_n) \ind_{\mX^\varepsilon_n}(\ux_n).
\end{equation}
Since their evolution is simply given by the global transport of $n$~particles, by which the equilibrium~$M_{\beta /2}^{\otimes n} \ind_{\mX^\varepsilon_n}$ is invariant, for all times~$t \geq 0$ the bound~(\ref{ineq:W_n0}) is propagated by the transport and remains true. Hence, taking the marginals we get
\begin{align*}
 W_n^{\varepsilon, (k)}(t, \uz_k, \ul_k) 
 & \leq \sum_{\ul^*_{n-k} \in \Lambda_{n-k}} \frac{( \lambda C_0)^{|\ul_{k}| + |\ul^*_{n-k}|} \mu^{n - |\ul_{k}| - |\ul^*_{n-k}|}}{\mZ_{\mu}} M_{\beta/2}^{\otimes k}(\uv_k) \left(\ind_{\mX^\varepsilon_n}\right)^{(k)}(\ux_k)\\
 &\leq (\lambda C_0 + \mu)^{n-k} \frac{(\lambda C_0)^{|\ul_{k}|} \mu^{k - |\ul_{k}|}}{\mZ_{\mu}}  M_{\beta/2}^{\otimes k}(\uv_k)  \ind_{\mX^\varepsilon_n}^{(k)}(\ux_k),
\end{align*}
where the factor~$(\lambda C_0 + \mu)^{n-k}$ stems from the sum over~$\ul^*_{n-k}$, using the binomial theorem. 
The main difference in the linear case, compared to the general non-linear one, is that the bound~\eqref{ineq:W_n0} on the canonical densities uses the invariant density~$M_{\beta/2}^{\otimes n}$, that passes to the $k$-th marginals to become~$M_{\beta/2}^{\otimes k}$, contrary to the constant $C_0^n$ in the general case.

Eventually, these bounds over the marginals of the canonical densities leaves the following a priori estimate for the correlation functions
\begin{align}
F_n^\varepsilon(t, \uz_n, \ul_n ) & \leq \frac{M_{\beta/2}^{\otimes n}}{\mu^{n - |\ul_n|} \lambda^{|\ul_n|} } \sum_{p\geqslant 0} \frac{1}{p!} \frac{(\lambda C_0 + \mu)^{p} (\lambda C_0)^{|\ul_{n}|} \mu^{n - |\ul_{n}|}}{\mZ_{\mu}}  \ind_{\mX^\varepsilon_{n+p}}^{(n)} \nonumber \\
& \leq \frac{ C_0^{|\ul_{n}|} M_{\beta/2}^{\otimes n} }{\mZ_{\mu}} \sum_{p\geqslant 0}  \frac{(\lambda C_0 + \mu)^{p} }{p!} \int \ind_{\mX^\varepsilon_{p}} \nonumber \\
& \leq \frac{ C_0^{|\ul_{n}|} M_{\beta/2}^{\otimes n} }{\mZ_{\mu}} \sum_{q,r\geqslant 0} \frac{(\lambda C_0)^q \mu^r }{q!r!} \int \ind_{\mX^\varepsilon_{q+r}} \label{eq:binomth}
\end{align}
where at line~\eqref{eq:binomth} we computed a direct binomial theorem. The key point to end our a priori estimate is now to control the remaining quotient implying the partition function~$\mZ_\mu$ and the slightly modified version of it, which is the following proposition.
\begin{prop} \label{prop:ZestimeeGC} There exists a constant $C_d$ depending only on the dimension such that for $\mu$ large enough, in our mixed Boltzmann-Grad scaling~\eqref{eq:scaling}, we have
\begin{equation*}
\frac{1}{\mZ_{\mu}}\sum_{q,r\geqslant 0} \frac{ (\lambda C_0)^q\mu^r}{q!r!} \int \ind_{\mX^\varepsilon_{q+r}} \leq C_d^{C_0 \lambda}. 
\end{equation*}
\end{prop}
The proof of this technical result is given in Section~\ref{app:partition_function}, using an explicit expansion of the partition function according to the cumulants of the exclusion.
Eventually, this leads to the following proposition, which is the main argument of our long time analysis.

\begin{prop}[A priori estimates for the correlation functions] \label{prop:aprioriestimates}
For any $n \in\N$ and $\varepsilon > 0$ in the mixed scaling~\eqref{eq:scaling}, one has
\begin{equation}
F_n^\varepsilon(t, \uz_n, \ul_n )  \leq C_0^{|\ul_{n}|}M_{\beta/2}^{\otimes n} \times C^{C_0 \lambda}. \label{ineq:apriori}
\end{equation}
\end{prop}

\subsection{Continuity estimates} 

Thanks to the a~priori estimates~(\ref{ineq:apriori}) given in the previous Section~\ref{sec:apriori}, all the correlation functions~$(F_n^\varepsilon)$ extended by~0 out of~$\mD_n^\varepsilon$ belong to the space~$\mF_{n,\beta/4}$ defined in~\eqref{def:foncspace}.

In the following, up to change the initial temperature~$\beta$ to~$4\beta$, to simplify the notation in the computation below we will assume that they belong to the space~$\mF_{n,\beta}$, with
\begin{equation} \label{ineq:Fwnorm}
\sup_{t \geqslant 0} \ \| F_n^\varepsilon(t) \|_{n, {\beta }} \leq  C_0^n C^{C_0 \lambda},
\end{equation}
thanks to the estimates above, and similarly for the initial proximity, from Proposition~\ref{prop:initial_proximity},
\begin{equation} \label{ineq:initial_proximity}
\left\| \ind_{\mX_n^\varepsilon}M_\beta^{\otimes n} \varphi_0^{\otimes \ul_n} - F_n^\varepsilon(0, \ul_n) \right\|_{n, \beta}\leq  C_0^n \varepsilon.
\end{equation}

Our whole long time derivation is allowed thanks to the fact that the a~priori bounds~\eqref{ineq:Fwnorm} are valid for every time with the same inverse temperature~$\beta$. Indeed, in the non-linear case, this parameter is downgraded over time until vanishing in finite time~\cite[Section~5]{2013newton}.

Technically, the downgrading of this norm parameter stems from the fact that to control the collision operators, the sum over the velocities is resorbed thanks to a fraction of the sub-Gaussian decreasing, eventually providing the following continuity estimates.
\begin{prop}[Continuity of the successive-collision operators] \label{prop:cont_est_Q}\ \\
There exists a constant $C_d$ depending only on the dimension such that for all~$n,s \in \N^*$, and all times~$t>0$, fixing tags~$\ul_{s} \in \Lambda_k$ and choosing two inverse temperature $\beta' < \beta$, we have
\begin{equation} \label{eq:contestQ}
\begin{array}{c}
f_{n+s} \in \mF_{n+s,\beta} \Rightarrow Q_{n, \ul_s}(t) f_{n+s} \in \mF_{n, \beta '} \textit{, with} \\ \ \\
\Bigl \| Q_{n, \ul_s}(t) f_{n+s} \Bigr \|_{n, \beta'} \leq e^{n} \left( \frac{C_d\ t}{\sqrt{\beta^d(\beta - \beta')}} \right)^s \|f_{n+s}\|_{n+s, \beta}.
\end{array}
\end{equation} 
\end{prop}
For fixed tags, the proof of this proposition is exactly the same as in the paper~\cite{fou24} about non-ideal Rayleigh gas, adapted from the article by Bodineau, Gallagher and Saint-Raymond, which derived the hydrodynamic limit of this gas~\cite{2016brownian}. Compared to this proof, we merely chose to write $\textstyle \beta' = \beta\left(1 - \frac{1}{b} \right)$ for better clarity in what follows. One can find it written in Foug\`{e}res' thesis~\cite{fouthese}.
 
\subsection{Adaptive tree pruning of the Dyson expansion for long times}  \label{sec:tree_pruning}

The continuity estimates presented in the previous section justify the wellposedness of the Dyson series~\eqref{eq:Dyson} for short times. To perform a derivation for long times, we will iterate the Dyson series, but each iteration will make appear a factor~$e^n$ stemming from the estimate of Proposition~\ref{prop:cont_est_Q}. These factors stack, so that, without further adjustment, at each iteration the successive times would decrease extremely fast and their sum would eventually be summable, leaving the derivation on a finite short time (see~\cite{fouthese} for details). The method we use here, introduced in~\cite{2016brownian}, consists in putting aside these stacking factors~$e^n$ and to control them by bounding the number of collisions appearing in the Dyson series.

More precisely, for a fixed time $t>1$ we will split the time interval $[0,t]$ into $K$ pieces and impose a piece-dependent amount of collisions on each small interval of this cutting. The number~$K$ will be tuned in the end of the proof, and henceforth we write
\begin{equation} \label{eq:time_split}
t = \sum_{i=1}^K h_i, \mbox{ with time steps } t^\mathrm{p}_k \doteq t - \sum_{j=1}^k h_i,
\end{equation}
so that $t_K^\mathrm{p} = 0$. Like in~\cite{fou24}, and contrary to \cite{2016brownian}, we will not choose a uniform cutting, but an adaptive one.
At the $k$-th time quantum of length $h_k$, we want at most $2^k$ collisions to have happened; we prune the collision tree every time it becomes more than exponentially big. Explicitly, between $t$ and $t_1^\mathrm{p} = t-h_1$,  we first truncate the Dyson series~(\ref{eq:Dyson}) to 2 collisions, then expand it again between $t-h_1$ and $t-h_2$ truncated to $2^2$ collisions, and iterate this process $K$ times.
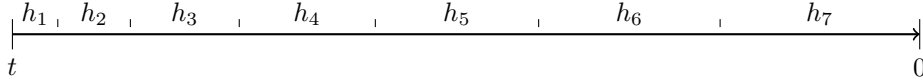
\begin{figure}[h!] 
\centering
\begin{tikzpicture}
\newcommand{\ho}{0.05}
\newcommand{\hd}{0.13}
\newcommand{\htr}{0.25}
\newcommand{\hc}{0.4}
\newcommand{\hp}{0.58}
\newcommand{\hs}{0.78}
\draw[black, thick, ->] (0,0) -- (12,0);
\draw[black] (12,0.2) -- (12,-0.2) node[below]{$0$};
\draw[black, thin, dashed] (\ho*12,0.2) -- (\ho*12,0);
\draw[black, thin, dashed] (\hd*12,0.2) -- (\hd*12,0) ;
\draw[black, thin, dashed] (\htr*12,0.2) -- (\htr*12,0);
\draw[black, thin, dashed] (\hc*12,0.2) -- (\hc*12,0);
\draw[black, thin, dashed] (\hp*12,0.2) -- (\hp*12,0);
\draw[black, thin, dashed] (\hs*12,0.2) -- (\hs*12,0);
\draw[black, thin, dashed] (\ho*6,0) node[above]{$h_1$};
\draw[black, thin, dashed] (\hd*6 + \ho*6,0) node[above]{$h_2$};
\draw[black, thin, dashed] (\htr*6+\hd*6,0) node[above]{$h_3$};
\draw[black, thin, dashed] (\hc*6 + \htr*6,0) node[above]{$h_4$};
\draw[black, thin, dashed] (\hp*6 + \hc*6,0) node[above]{$h_5$};
\draw[black, thin, dashed] (\hs*6 + \hp*6 ,0) node[above]{$h_6$};
\draw[black, thin, dashed] (\hs*6 + 6 ,0) node[above]{$h_7$};
\draw[black] (0, 0.2) -- (0,-0.2) node[below]{$t$};
\end{tikzpicture}
\caption{Backwards division of the time interval under study}
\label{fig4}
\end{figure}

\noindent
We denote the step number of added tagged particles and the step total number of particles 
\begin{equation} \label{def:J_notation}
L_k \doteq \sum_{i=1}^k |\ul^*_{j_i}| \ \ \und\ \  N_k \doteq n + j_1 + \cdots + j_k.
\end{equation}
This yields, similarly as in~\cite{2016brownian, fou24}, the following pruned expansion
\begin{align} \label{eq:pruning}
F_n^\varepsilon(t) & = \sum_{\ \ \left(\substack{j_i \leqslant 2^i \\ \ul^*_{j_{i}} \in \Lambda_{j_i} } \right)_{1 \leqslant i \leqslant K}}  p_\mu^{L_K}  Q_{n, \ul^*_{j_1}}(h_1) \dots  Q_{N_{K-1}, \ul^*_{j_{K}}}(h_{K}) F_{N_K}^\varepsilon(0) + R_n^{[K]}(t)
\end{align}
where the remainder is defined as
\begin{equation} \label{eq:remK}
R_n^{[K]}(t) \doteq \sum_{k=1}^K \sum_{\left(\substack{j_i \leqslant 2^i \\ \ul^*_{j_{i}} \in \Lambda_{j_i} } \right)_{i \leqslant k-1}}  Q_{n, \ul^*_{j_1}}(h_1) \dots  Q_{N_{k-2}, \ul^*_{j_{k-1}}}(h_{k-1})  \sum_{j_k>2^k} \sum_{\ul^*_{j_k} \in \Lambda_{j_k}} p_\mu^{L_k} Q_{N_{k-1}, \ul^*_{j_k}}(h_{k}) F_{N_{k}}^\varepsilon(t^{\mathrm{p}}_k).
\end{equation}
Based on the limit expansion~\eqref{eq:Dyson_lim}, one can write the same decomposition for the limit family~$(G_n)_{n\geqslant 0}$:
% $\overset{ _{\lim}}{Q} 
\begin{align} \label{eq:pruning_limit}
G_n(t) & =\sum_{\ \ \left(j_i \leqslant 2^i \right)_{1 \leqslant i \leqslant K}}   Q^{\lim}_{n, \underline{0}_{j_1}}(h_1) \dots  Q^{\lim}_{N_{K-1}, \underline{0}_{j_{K}}}(h_{K}) G_{N_K}(0) + R_n^{[K], \lim}(t).
\end{align}
Let us denote~$\widehat{G}_n(t) \doteq G_n(t) - R^{[K], \lim}_n(t)$ the limit pruned expansion.
Since the chosen condition of a sub-exponential number of collisions is very restrictive at first, and then gradually relaxed, the adaptive cutting times are chosen small at first and then progressively bigger (see Fig.~\ref{fig4}). The bound on the pruned-out remainder is given in the following proposition.

\begin{prop}[Estimate of the pruned-out term] \label{prop:estrem}
With the previous notation, for any choice of power~$\alpha \in (0,3/4)$, and $K$ large enough satisfying $t \lesssim K^{\frac{3}{4}-\alpha}$, a good choice of time cutting \[\uh = (h_1, \dots, h_K)\] provides the following estimate 
\begin{equation} \label{ineq:est_pruned_out}
\lnorm R_n^{[K]}(t) \rnorm_{\Lp^\infty(\mD^d)} + \lnorm R_n^{[K], \lim}(t) \rnorm_{\Lp^\infty(\mD^d)}  \leq C^{C_0 \lambda} n^{cn} e^{-2^{K - K^\alpha}}. 
\end{equation}
\end{prop} 
 
Note that this estimate imposes a technical condition on the mixed scaling of~$\lambda$ and $K$ for the error to be small; we will tune this scaling in Section~\ref{sec:conclusion1}.
The factor $n^{cn}$ stems from the fact that we generalize the result to all the correlation functions, and not only the first one.  Indeed, one will see that this factor is due to the usual bound $C^n$ stacking at each iteration of the cutting.

\begin{prof}
The proof is very similar to the one found in~\cite{fou24}. We give it for the hard~sphere version, since the limit version is identical. Using the a~priori estimates~(\ref{ineq:apriori}) on the densities and the continuity estimate on successive-collision operators given in Proposition~\ref{prop:cont_est_Q}, at given $(\ul^*_{j_i})_{i \leq k-1}$ one has for every~$k\in \lbr 1, K \rbr$
\begin{align*}
\lnorm \sum_{j_k>2^k} \sum_{\ul^*_{j_k}\in \Lambda_{j_k}} p_\mu^{L_k} Q_{N_{k-1}, \ul^*_{j_k}}(h_{k}) F_{N_{k}}^\varepsilon(t^{\mathrm{p}}_k) \rnorm_{N_{k-1}, \beta /2} 
& \leq e^{N_{k-1}} \sum_{j_k>2^k} \sum_{\ul^*_{j_k}\in \Lambda_{j_k}} p_\mu^{L_k} \left( \frac{ \sqrt{2} C_d  h_k}{\beta^{(d+1)/2}} \right)^{j_k} \lnorm  F^\varepsilon_{N_{k} }(t^\mathrm{p}_k)\rnorm_{N_{k}, \beta}\\
& \leq e^{N_{k-1}} \sum_{j_k>2^k} \sum_{\ul^*_{j_k}\in \Lambda_{j_k}} p_\mu^{L_k}  \left( \frac{ \sqrt{2} C_d  h_k}{\beta^{(d+1)/2}} \right)^{j_k} C^{N_{k}} C^{C_0 \lambda}.
\end{align*}
For $\mu$ large enough in the scaling~\eqref{eq:scaling}, we have $p_\mu \leq 1$, so that the sum over $\ul^*_{j_k}\in \Lambda_{j_k}$ only gives a factor $2^{N_k}$ that can be resorbed in the term $C^{N_k}$. We then iterate Proposition~\ref{prop:cont_est_Q}, downgrading the parameter $\beta/2$ by $\beta/(4k)$ at each step, so that it remains greater than~$\beta/4$. Hence we can write, grouping~$C^{N_{k-1}}$ and all the appearing terms of the form $e^{N_i}$ together as a power of a constant $\hat{C}$, that
\begin{align}
& \lnorm Q_{n, \ul^*_{j_1}}(h_1) \dots  Q_{N_{k-2}, \ul^*_{j_{k-1}}}(h_{k-1})  \sum_{j_k>2^k} \sum_{\ul^*_{j_k} \in \Lambda_{j_k}} p_\mu^{L_k} Q_{N_{k-1}, \ul^*_{j_k}}(h_{k}) F_{N_{k}}(t^{\mathrm{p}}_k) \rnorm_{n, \beta/4} \label{eq:succbound} \\  
& \hspace*{3mm} \leq C^{  C_0 \lambda } \hat{C}^{\sum_{i=0}^{k-1} N_i} \left( \left(\frac{4}{\beta}\right)^\frac{d+1}{2} \sqrt{4k} C_d  h_1 \right)^{j_1} \dots  \left( \left(\frac{4}{\beta}\right)^\frac{d+1}{2} \sqrt{4k} C_d  h_{k-1}  \right)^{j_{k-1}} \sum_{j_k>2^k} \left( \frac{\sqrt{2} C_d  h_k}{\sqrt{\beta}} \right)^{j_k} \cdotp \nonumber
\end{align}
We now observe, on the one hand, that recalling notation~\eqref{def:J_notation} for $N_i$ and since for $i \leq k-1$, $j_i \leq 2^i$, we have $\textstyle \sum_{i=0}^{k-1} N_i \leq nk + 2^{k+1}$. On the other hand, we can also put aside from the sum the following factors
\begin{align*}
\left(\left(\frac{4}{\beta}\right)^\frac{d+1}{2} 2 k^{\frac{1}{4} - \alpha} C_d  \right)^{\sum_{i=0}^{k-1} j_i} & \leq (\hat{C}_\beta k^{\frac{1}{4} - \alpha})^{2^{k}}.
\end{align*}
This computation is what allows us to gain a power $\textstyle \frac{1}{4}$ on the time scaling that we imposed in Theorem~\ref{theo1}, compared to~\cite{fou24}. Indeed, in the following, the computation above harnesses the full power decay of the last time interval~$h_k$, whereas some of it was lost in~\cite{fou24}. In the end, for a possibly larger constant~$C$ depending on~$d$ and~$\beta$, we~get
\begin{align*} \label{absorbtion} 
\lnorm R_n^{[K]}(t) \rnorm_{\Lp^\infty} & \leq C^{C_0 \lambda} \sum_{k=1}^K C^{nk} (Ck^{\frac{1}{4}-\alpha})^{2^k} \sum_{j_1 = 0}^{2} \cdots \sum_{j_{k-1} = 0}^{2^{k-1}}   (k^{\frac{1}{4}+\alpha} h_1)^{j_1} \dots   (k^{\frac{1}{4}+\alpha} h_{k-1})^{j_{k-1}} \sum_{j_k>2^k} \left(C  h_k \right)^{j_k}.
\end{align*}
Eventually, we consider similarly as in~\cite{fou24}, for all $1 \leq i \leq K$,
\begin{equation} \label{def:htilde}
\tilde{h}_i \doteq  \frac{e^{-2^{(K-K^{1- \alpha } - i)}}}{2 K^{\frac{1}{4}+\alpha}} \leq \frac{1}{2K^{\frac{1}{4}+\alpha}},
\end{equation}
and renormalize them such that
\begin{equation*}
h_i \doteq \frac{t}{\sum_{j=1}^K \tilde{h}_j} \tilde{h}_i \leq \tilde{h}_i.
\end{equation*}
Indeed, as soon as~$t \lesssim K^{\frac{3}{4} - 2\alpha}$, we can write
\begin{align*}
\frac{1}{t}\sum_{i = 1}^K \tilde{h}_i 
& \geq \frac{1}{t} \sum_{j = 0}^{\lfloor K^{1-\alpha} \rfloor} \tilde{h}_{K-j} \\
& \geq \frac{1}{t} \sum_{j = 0}^{\lfloor K^{1-\alpha} \rfloor}\frac{e^{-1}}{2 K^{\frac{1}{4}+\alpha} }  \geq 1.
\end{align*}
Now, the time interval lengths~$(h_i)_{1\leqslant i \leqslant K}$ cover~$t$ as imposed by~\eqref{eq:time_split}, and their choice provides, summing the geometric series over~$(j_i)_{1\leqslant i \leqslant K}$,
\begin{align*}
\lnorm R_n^{[K]}(t) \rnorm_{\Lp^\infty} &  \leq C^{C_0\lambda} \sum_{k=1}^K C^{nk} (Ck^{\frac{1}{4}-\alpha})^{2^k} \times 2^k \sum_{j_k>2^k}  \left( \frac{ e^{-2^{K-K^{1-\alpha}-k}}}{2 K^{\frac{1}{4}+\alpha}} \right)^{j_k}\\
& \leq C^{C_0\lambda} \sum_{k=1}^K C^{nk} \left( \frac{CK^{\frac{1}{4}-\alpha}}{K^{\frac{1}{4}+\alpha}}\right)^{2^k} \left( \frac{e^{-2^{K-K^{1-\alpha}-k}}}{2}\right)^{2^k}.
\end{align*}
Observe now that there exists a constant $c$ such that $(k \geq c \log n) \Rightarrow (nk \leq 2^k)$, so that in this case the factor $C^{nk}$ is absorbed by $C^{2^k}$, and for $k \leq c \log n$, then $C^{nk} \leq C^{cn \log n} = n^{(c \log C) n}$. Now, the denominator $K^{\frac{1}{4}+\alpha}$ crushes the term $C K^{\frac{1}{4} - \alpha}$ for $K$ large enough and we end up with
\begin{align*}
\lnorm R_n^{[K]}(t) \rnorm_{\Lp^\infty} \leq C^{C_0\lambda} \tilde{C}^{n \log n}   e^{-2^{K-K^{1-\alpha}}}. 
\end{align*}
This completes the proof of Proposition~\ref{prop:estrem}.
 \CQFD
\end{prof}

\subsection{Discarding trajectories implying several labelled particles} 
 \label{sec:disc_labelled}
Since the limit pseudo-trajectories defined from the limit hierarchy~\eqref{eq:pseudo_traj_form_lim} only imply collisions with particles at equilibrium (which are in wide majority), we must get rid of the ones including collisions with perturbed tagged particles. We hence write our pruned expansion as
\begin{align}
& \sum_{\ \ \left(\substack{j_i \leqslant 2^i \\ \ul^*_{j_{i}} \in \Lambda_{j_i} } \right)_{1 \leqslant i \leqslant K}} p_\mu^{L_K}  Q_{n, \ul^*_{j_1}}(h_1) \dots  Q_{N_{K-1}, \ul^*_{j_{K}}}(h_{K}) F^\varepsilon_{N_K}(0) \nonumber \\
&  \hspace{14mm} = \sum_{\left(j_i \leqslant 2^i \right)_{1 \leqslant i \leqslant K}} Q_{n, \underline{0}_{j_1}}(h_1) \dots  Q_{N_{K-1}, \underline{0}_{j_{K}}}(h_{K}) F^\varepsilon_{N_K}(0) \label{def:Fmain}\\
& \hspace{18mm} + \sum_{ \left(j_i \leqslant 2^i \right)} \sum_{k=1}^K \sum_{\substack{(\ul^*_{j_i})_{i \leqslant K} \\ \ul^*_{j_{k}} \neq \underline{0}_{j_k}}} p_\mu^{L_K}  Q_{n, \ul^*_{j_1}}(h_1) \dots  Q_{N_{K-1}, \ul^*_{j_{K}}}(h_{K}) F^\varepsilon_{N_K}(0), \label{def:F_UE} \\
&  \hspace{18mm} \doteq \widehat{F}_n^\varepsilon(t) + F^{\varepsilon, \mathrm{u.e.}}_n(t), \nonumber
\end{align}
where we denote~$\widehat{F}_n^\varepsilon(t)$ the main term~\eqref{def:Fmain} containing only collisions with equilibrium, and $F^{\varepsilon, \mathrm{u.e.}}_n(t)$ the one implying unwanted encounters~\eqref{def:F_UE}. We bound the latter in the following proposition.
\begin{prop}[Encountering tagged particles is rare] \label{prop:discard} 
As long as \begin{equation*}
p_\mu \leq \frac{1}{2^K},
\end{equation*}
in the same time setting as in Proposition~\ref{prop:estrem}, the unwanted pseudotrajectories implying tagged particles are bounded by
\begin{equation} \label{ineq:discard_tagged}
\| F^{\varepsilon, \mathrm{u.e.}}_n(t)\|_{\Lp^\infty(\mD^d)} \leq  C^{nK + A^K}\ p_\mu.
\end{equation}
\end{prop}
The reader may think of the factor~$C^{A^K}$ as a small negative power of~$\varepsilon$ in the final scaling, which we will compute in the following Section~\ref{sec:conclusion1}: the scaling proportion~$p_\mu$ will have to compensate it.

\begin{propf} A computation directly adapted from the initial proximity bound~\eqref{ineq:initial_proximity} leads to the estimate
\begin{equation*}
\left\| F_{N_K}^\varepsilon(0, \tilde{\ul}_{N_K}) \right\|_{N_K, \beta} \leq  (C_0 )^{N_K},
\end{equation*}
so that, using the binomial identity
\begin{equation*}
\sum_{\ul^*_{j_i} \in \Lambda_{j_i}} p_\mu^{|\ul^*_{j_i}|} = (1+p_\mu)^{j_i}
\end{equation*}
and the same continuity estimates on the successive-collision operators as in the proof of Proposition~\ref{prop:estrem}, one has
\begin{align*}
\| F^{\varepsilon, \mathrm{u.e.}}_n(t)\|_{\Lp^\infty(\mD^d)} & \leq  C^{nK} (CK^{\frac{1}{4}-\alpha})^{2^K} \sum_{(j_i \leqslant 2^i)_{i \leqslant K}} \sum_{k=1}^K \left( (1 + p_\mu)^{j_k}  - 1 \right)  \prod_{i=1}^K \left(\frac{1}{2} \right)^{j_i} .
\end{align*}
Notice on the one hand that 
\begin{align*}
(CK^{\frac{1}{4}-\alpha})^{2^K} & = \exp\left[ 2^K \log(CK^{\frac{1}{4}-\alpha})\right] \\
& \leq \exp\left(A^K \right)
\end{align*}
for any $A > 2$ as long as $K$ is large enough (depending on $A$). On the other hand, using $j_k \leq 2^k$ and taking $\textstyle p_\mu \leq 2^{-K} \leq 2^{-k}$, we have by convexity on~$[0, 2^{-k}]$ that
\begin{align*}
(1 +  p_\mu)^{j_k}  - 1 & \leq (1 + p_\mu)^{2^k}  - 1 \\
& \leq (e-1) 2^k  p_\mu,
\end{align*}
so that eventually
\begin{align*}
\| F^{\varepsilon, \mathrm{u.e.}}_n(t)\|_{\Lp^\infty(\mD^d)} & \leq   (e-1) p_\mu \times   C^{nK + A^K}  \sum_{k=1}^K 2^k \sum_{(j_i \leqslant 2^i)_{i \leqslant K}}  \prod_{i=1}^K \left(\frac{1}{2} \right)^{j_i}\\
&\leq p_\mu\ \widehat{C}^{nK + A^K}
\end{align*} 
for another constant $\widehat{C}$ absorbing the factor $(e-1) 2^{K+1} \times 2^K$, concluding the proof. \CQFD
\end{propf}

\subsection{Proof of the convergence} \label{sec:conclusion1}
Now that the pruned-out terms have been controlled, we want to compare the pruned terms~$\widehat{F}^\varepsilon_n(t)$ and~$\widehat{G}_n(t)$, which is the last step before our choice of scaling for $K$ and the conclusion of Theorem~\ref{theo1}. As explained in Section~\ref{sec:strategy}, our strategy relies on considering pseudo-trajectories without recollisions.
This method is an adaptation of~\cite[Section~5]{2016brownian} which is now classical; it follows from several approximations: an energy truncation and a time separation of the collisions are operated, so as to be able to construct a small set of bad collision parameters, outside of which there will be no recollision.

This time separation method does not yield the best quantitative estimates, but in our study in long time the worst error is made with the pruned-out terms, so that here we use this method anyway. This allows concision on the one hand, and on the other hand one can thus compare it with the optimized method presented in the companion paper~\cite{fou26cumulants}.

Like in~\cite{fou24}, we still refine the quantitative bounds using the same control on\[Q_{n, \ul^*_{j_1}}(h_1) \dots  Q_{N_{K-1}, \ul^*_{j_{K}}}(h_{K})\] as in Propositions~\ref{prop:estrem} and~\ref{prop:discard}, which induces our factor $ {C_\beta}^{nK + A^K} $, instead of the original crude bound with $|Q|_{1, J_{k}}(t)$, which gave a factor $(C t)^{2^K}$ (see~\cite{2016brownian}). From a physical perspective, we decompose the time interval into small pieces whose lengths are adapted to the maximum number of particles that may appear in them, so that the dynamics behaves similarly during each one of them. Hence, as long as time does not get too big with respect to the number of pieces, none of the estimates depends on the total time length.

\paragraph{Pseudo-trajectory formulation}
First of all, let us notice that the pruned expansion~\eqref{eq:pruning} has a pseudo-trajectory formulation similar to the original one~\eqref{eq:pseudo_traj_form}, summing over the successive numbers of collisions~$(j_i \leq 2^i)_{i\leq K}$, with the following additional condition on the collision times, located between the time steps~\eqref{eq:time_split}:
\begin{equation*}
\underline{t}_{J_K} \in T_{\uj_{K} }^{}(t) \doteq \Bigl \{ (t_1, \dots, t_{J_K})\in T_{J_K}(t)\ ,\  \{ t_{J_k +1}, \dots, t_{J_{k+1}} \} \subset [t_{k+1}^\mathrm{p}, t_{k}^\mathrm{p}] \Bigr \},
\end{equation*}
where $J_K \doteq N_K - n = j_1 + \dots + j_K$ denotes the step number of collisions. As announced above, we start by computing a few approximations, initiated by an energy truncation, so as to work with bounded velocities.

\paragraph{Energy truncation} 
We hence consider the following functionals with truncated energy
\begin{align}
\label{eq:energy_truncation}
\widehat{F}_{n}^{\varepsilon, [\bV]}(t)  \doteq \sum_{\uj_K} \sum_{\substack{\uchi_{J_K} \\ L_K = 0 }}  \int_{T^{}_{\underline{j}_K}(t)} \dd \ut_{J_K} \int \dd \underline{\omega}_{J_K} \dd v_{n+1} \dots \dd v_{N_K} & \prod_{i=1}^{{J_K}} s_i \langle \omega_{i},  v_{n+i} - v_{m_i}^{[t_i^+]} \rangle_+ \\ 
& \times F_{N_K}^\varepsilon \bigl(0, \uz_{N_K}^{[0]} \bigr) \ind_{\left\|\uv_{N_K}^{[0]}\right\|^2 \leqslant \bV^2}, \nonumber
\end{align} 
and similarly the truncated limit functions $\widehat{G}_n^{[\bV]}$. The error made by truncating the velocities is
\begin{align*}
\widehat{F}_{n}^{\varepsilon}(t) - \widehat{F}_{n}^{\varepsilon, [\bV]}(t) = \sum_{ \left(j_i \leqslant 2^i \right)} Q_{n, \underline{0}_{j_1}}(h_1) \dots  Q_{N_{K-1}, \underline{0}_{j_{K}}}(h_{K}) \left[ F^\varepsilon_{N_K}(0) \ind_{\left\|\uv_{N_K}^{[0]}\right\|^2 \geqslant \bV^2} \right],
\end{align*}
with the following estimate on the initial functional on the right:
\begin{align*}
\left\| F_{N_K}^\varepsilon \bigl(0, \uz_{N_K}^{[0]}, \tilde{\ul}_{N_K} \bigr) \ind_{\left\|\uv_{N_K}^{[0]}\right\|^2 \geqslant \bV^2} \right\|_{N_K, \beta/2} &  \leq \sup_{ \uu_{N_K}\in \R^{d N_K}} \left\lvert F_{N_K}^\varepsilon \bigl(0 \bigr) e^{\beta \left\|\uu_{N_K}\right\|^2 } e^{-\frac{\beta}{2}\left\|\uu_{N_K}\right\|^2} \ind_{\left\|\uu_{N_K}\right\|^2 \geqslant \bV^2} \right\rvert \nonumber \\
& \leq \lnorm F_{N_K}^\varepsilon(0) \rnorm_{N_K, \beta } e^{-\frac{\beta}{2} \bV^2}.
\end{align*}
Hence, applying the same bounds as in Propositions~\ref{prop:estrem} and~\ref{prop:discard} in the same setting, we end up with the following lemma.
\begin{lemm}[Energy truncation error] The error due to the energy truncation is bounded by
\begin{align} \label{ineq:energy_error}
\|\widehat{F}_{n}^{\varepsilon}(t) - \widehat{F}_{n}^{\varepsilon, [\bV]}(t) \|_{\Lp^\infty} & \leq C^{nK + A^K} \exp\left(-\frac{\beta}{2} \bV^2 \right).
\end{align}
The same holds for its limit version~$\widehat{G}_{n}^{}(t) - \widehat{G}_{n}^{ [\bV]}(t)$.
\end{lemm}

\paragraph{Time separation}
We now need the successive collisions to be separated enough in time, to avoid pathological geometric recollisions. Like in the previous section, let us define
\begin{align}\label{eq:time_sep}
\widehat{F}_{n}^{\varepsilon, [\bV, \delta]}(t)  \doteq \sum_{\uj_K} \sum_{\substack{\uchi_{J_K} \\ L_K = 0 }} \int_{T^{[\delta]}_{\underline{j}_K}(t)} \dd \ut_{J_K} \int \dd \underline{\omega}_{J_K} \dd v_{n+1} \dots \dd v_{N_K} & \prod_{i=1}^{{J_K}} s_i \langle \omega_{i},  v_{n+i} - v_{m_i}^{[t_i^+]}\rangle_+ \\ 
& \times F_{N_K}^\varepsilon \bigl(0, \uz_{N_K}^{[0]} \bigr) \ind_{\left\|\uv_{N_K}^{[0]}\right\|^2 \leqslant \bV^2}, \nonumber
\end{align}
with the separation condition encoded in the following time set, over which we integrate,
\begin{equation} \label{eq:Time_sep}
T^{[\delta]}_{\underline{j}_K}(t) = \Bigl \{ \ut \in T^{}_{\underline{j}_K}(t) \ , \ t_{i-1} - t_i > \delta \Bigr \}.
\end{equation}
The limit version~$\widehat{G}_n^{[\bV, \delta]}$ is defined by the same time restriction. The error of time separation is
\begin{align*}
\widehat{F}_n^{\varepsilon, [\bV]} - \widehat{F}_n^{\varepsilon, [\bV, \delta]} \doteq \sum_{\uj_K} \sum_{\substack{\uchi_{J_K} \\ L_K = 0 }} & \int_{\left(T^{[\delta]}_{\underline{j}_K}(t)\right)^c} \dd \ut_{J_K} \int   \dd \underline{\omega}_{J_K} \dd v_{n+1} \dots \dd v_{N_K} \\
& \times  \prod_{i=1}^{{J_K}} s_i \langle \omega_{i}, v_{m_i}^{[t_i^+]} - v_{n+i} \rangle_+ 
 \times F_{N_K}^\varepsilon \bigl(0, \uz_{N_K}^{[0]} \bigr) \ind_{\left\|\uv_{N_K}^{[0]}\right\|^2 \leqslant \bV^2}, \nonumber
\end{align*} 
where one can write 
\begin{align} \label{eq:union_comp}
\left(T^{[\delta]}_{\underline{j}_K}(t)\right)^c = \bigcup_{i = 1}^{J_K-1} \Bigl \{ \ut \in T^{}_{\underline{j}_K}(t) \ , \ t_{i} - t_{i+1} \leqslant \delta \Bigr \}.
\end{align}
Now, the integral in time over one of these sets, using the same method as before, changes the estimate of the corresponding successive-collision operator from~$\frac{h_k^{j_k}}{j_k!}$ to~$\frac{\delta h_k^{j_k - 1}}{(j_k - 1)!}$. Since the loss of~$\textstyle \frac{h_k}{j_k}$ and the factor~$(J_K - 1 )$ coming from the union~\eqref{eq:union_comp} easily resorb into the bigger factor~$C^{A^K}$, we get in the end the following lemma.
\begin{lemm}[Time separation error] For this error, we have the following estimate
\begin{align} \label{ineq:time_separation}
\|\widehat{F}_n^{\varepsilon, [\bV]} - \widehat{F}_n^{\varepsilon, [\bV, \delta]} \|_{\Lp^\infty} & \leq C^{nK + A^K} \delta.
\end{align}
 Once again, the same holds for its limit version~$\widehat{G}_n^{ [\bV]} - \widehat{G}_n^{ [\bV, \delta]}$ for similar reasons.
\end{lemm}

\paragraph{Restriction to non-pathological collision parameters}
Finally, at fixed~$\varepsilon$, we have to restrict the collision parameters to non-pathological configurations, leading to no recollision during the transport flow. First, to compute the convergence of the $n$-th marginal, we have to consider final configurations $\uz_n = \uz_n^{[t]}$ that do not directly lead to  recollisions, i.e. that belong to the following set of past-excluding configurations, with some extra room $\varepsilon_d \doteq \varepsilon^\frac{d}{d+1}$:
\begin{equation} \label{def:past-exclu}
\mE^t_n(\varepsilon_d) \doteq \Bigl \{ \uz_n \in \mD^\varepsilon_n\ \Bigl \lvert \ \forall \tau \in [0,t], \ (\ux_n - \tau \uv_n) \in \mD^{\varepsilon_d}_n \Bigr \}.
\end{equation}
This set is the whole domain for $n=1$, like in~\cite{2016brownian, fou24}.
Then, the restriction on the collision parameters is done based on the following geometric result, proved in~\cite{2013newton, 2016brownian} and formalizing the arguments of Lanford~\cite{1975lanford}, based on billiards theory.

Similarly to the notation~$\uz_n^{[\tau]} = (\ux_n^{[\tau]}, \uv_n^{[\tau]})$ for the hard sphere pseudo-trajectories, let us denote~$\uze_n^{[\tau]} = (\uy_n^{[\tau]}, \uu_n^{[\tau]})$ the limit pseudo-trajectories, and $N[\tau]$ the number of particles in the trajectory at time $\tau \in [0, t]$, from $N[t] = n$ to $N[0] = N_K$. Hence, the following lemma asserts that once the previous truncations computed, choosing collision parameters away from a set of small measure, the hard sphere and limit pseudo-trajectories are easy to compare.
\begin{lemm} \label{lemm:geom}
Considering a history~$(\uj_K, \uchi_{J_K})$ with collision times $\ut_{J_K} \in T^{[\delta]}_{\uj_K}(t)$ ($\delta$-separated) and a final configuration~$\uz_n^{[t]}$, given a maximum energy~$\bV^2 > 0$, there exists a set of pathological collision parameters
\begin{equation*}
\Pi(\uz_n^{[t]}, \uj_K, \uchi_{J_K}) \subset ( \Sf^{d-1} \times \R^d )^{J_K}
\end{equation*} 
with small volume
\begin{equation} \label{ineq:volume_Pi}
|\Pi(\uz_n^{[t]}, \uj_K, \uchi_{J_K})| \leq C J_K N_K \left( \varepsilon^{\frac{d}{d+2}} + \bV^d \times \varepsilon^{\frac{d-1}{3(d+1)}} + \bV^\frac{d+1}{2} \left( \frac{\varepsilon^{\frac{d}{d+1}}}{\delta} \right)^\frac{d-1}{2} \right), 
\end{equation}
and such that, assuming
\begin{enumerate}[label=\roman*.]
\item the collision parameters are non-pathological: $(\uom_{J_K}, v_{n+1}, \dots, v_{N_K}) \notin \Pi(\uz_n^{[t]}, \uj_K, \uchi_{J_K})$
\item the energy of the corresponding pseudo-trajectory remains bounded: $\|\uv_{N_K}^{[0]}\|^2 \leq \bV^2$ 
\item the final configuration is past-excluding for the free-flow: $ \uz_n^{[t]} \in \mE^t_n(\varepsilon_d)$,
\end{enumerate}
then
\begin{enumerate}
\item the hard sphere pseudo-positions remain sufficiently far away: $\forall \tau \in [0,t], \ \ux_{N[\tau] }^{[\tau]} \in \mD^{\varepsilon_d / 2}_{N[\tau]}$
\item the \emph{velocities} of the hard sphere and limit trajectories coincide:  $\forall \tau \in [0,t], \ \uv_{N[\tau] }^{[\tau]} = \uu_{N[\tau]}^{[\tau]}$
\item the \emph{positions} of both trajectories remain close: $ \forall \tau \in [0,t], \forall i \leq N[\tau], d(x_i^{[\tau]}, y_i^{[\tau]}) \leq J_K \varepsilon$.
\end{enumerate}
\end{lemm}
Let us observe that~\emph{3.} is a consequence of~\emph{2.} since, when the velocities coincide, the only difference between the positions is the shift of~$\varepsilon \omega_i$ that happens at each particle adjunction. Furthermore, \emph{2.} is a consequence of~\emph{1.} since when the collision parameters are the same, and if the particles do not collide between the particle adjunctions, the velocities of both pseudo-trajectories are identically determined. 
This result has been proved~\cite{2013newton, 2016brownian} for a one-species gas, but the dynamics is strictly identical for our mixture gas.

Hence, we consider the following functional restricted to non-pathological collision parameters
\begin{align}\label{eq:path_traj}
\widetilde{F}_{n}^{\varepsilon, [\bV, \delta]}(t)  \doteq \sum_{\uj_K} \sum_{\substack{\uchi_{J_K} \\ L_K = 0 }}  & \int_{T^{[\delta]}_{\underline{j}_K}(t)} \dd \ut_{J_K}  \int_{ \Pi(\uz_n^{[t]}, \uj_K, \uchi_{J_K})^c} \dd \underline{\omega}_{J_K} \dd v_{n+1} \dots \dd v_{N_K} \\ 
& \times \prod_{i=1}^{{J_K}} s_i \langle \omega_{i},  v_{n+i} - v_{m_i}^{[t_i^+]} \rangle_+  \times F_{N_K}^\varepsilon \bigl(0, \uz_{N_K}^{[0]} \bigr) \ind_{\left\|\uv_{N_K}^{[0]}\right\|^2 \leqslant \bV^2}, \nonumber
\end{align}
and its limit version~$\widetilde{G}_n^{[\bV, \delta]}$. Now, the error 
$\widetilde{F}_n^{\varepsilon, [\bV, \delta]} - \widehat{F}_n^{\varepsilon, [\bV, \delta]}$ is supported on $\Pi(\uz_n^{[t]}, \uj_K, \uchi_{J_K})$, so that we will use the control on its volume~\eqref{ineq:volume_Pi} to control the successive-collision operators, concluding the bounds with the usual computation. More precisely, in the proof~\cite{fou24} of Proposition~\ref{prop:cont_est_Q}, we bound the collision operators~\eqref{def:collision_op2} in the following way:  	
\begin{align*}
e^{\beta \|\uv_j\|^2} \Bigl| \mC_{j}^\ell f_{j+1}  \Bigr| & \leq \sum_{i=1}^j \int \dd \omega_j \dd v_{j+1} \left| v_{j+1} - v_i \right| \cdot \|f_{j+1} \|_{j+1, \beta'} \times e^{-(\beta'- \beta) \|\uv_{j}\|^2 - \beta' |v_{j+1}|^2 } \\
& \leq \|f_{j+1} \|_{j+1, \beta'} \int \dd \omega_j \dd v_{j+1} \left( j \left| v_{j+1} \right| + \sum_{i=1}^j \left|v_i\right| \right)  e^{-(\beta'- \beta) \|\uv_{j}\|^2 - \beta' |v_{j+1}|^2 }  .
\end{align*}
By the Cauchy-Schwarz inequality, we used to bound the quantity above as
\begin{align*}
e^{\beta \|\uv_j\|^2} \Bigl| \mC_{j}^\ell f_{j+1}  \Bigr| & \leq \|f_{j+1} \|_{j+1, \beta'} \int \dd \omega_j \dd v_{j+1} \left( j \left| v_{j+1} \right| + \sqrt{ \frac{j}{2e (\beta'-\beta)} } \right)  e^{-\beta' |v_{j+1}|^2} \\
& \leq  \|f_{j+1} \|_{j+1, \beta'} \left( j \frac{c_d}{\sqrt{{\beta'}^{d+1}}} + \sqrt{ \frac{j}{2e (\beta'-\beta)} } \frac{\tilde{c}_d}{\sqrt{{\beta'}^{d}}} \right),
\end{align*}
for some constants~$c_d, \tilde{c}_d$ depending only on the dimension. For the present bound, one can write instead
\begin{equation*}
\int \dd \omega_{j} \dd v_{j+1} \left(j|v_{j+1}| + \sqrt{ \frac{j}{2e(\beta' - \beta)}} \right) e^{ - \beta' |v_{j+1}|^2} \leq \int \dd \omega_{j} \dd v_{j} \left( \frac{1}{\sqrt{2e\beta}} +  \sqrt{\frac{e^{-1}}{\beta - \beta'}} \right),
\end{equation*}
making appear the volume of collision parameters bounded in~\eqref{ineq:volume_Pi}, leading to the following estimate (the factors~$J_K$ and~$N_K$ resorb in the bigger factor~$C^{nK + A^K}$, for a slightly different constant~$C$ depending on~$\beta$).
\begin{lemm}[Error of restriction to non-pathological collision parameters] The restriction of considering only collision parameters chosen so as to avoid recollisions, leads to an error of order
\begin{align} \label{ineq:path_traj}
\|\widetilde{F}_n^{\varepsilon, [\bV, \delta]} - \widehat{F}_n^{[\bV, \delta]} \|_{\Lp^\infty} & \leq C^{nK + A^K} \left( \varepsilon^{\frac{d}{d+2}} + \bV^d \times \varepsilon^{\frac{d-1}{3(d+1)}} + \bV^\frac{d+1}{2} \left( \frac{\varepsilon^{\frac{d}{d+1}}}{\delta} \right)^\frac{d-1}{2} \right),
\end{align}
and similarly for the limit~$\widehat{G}_n^{ [\bV]} - \widehat{G}_n^{ [\bV, \delta]}$ from the same computation.
\end{lemm}

\paragraph{Harnessing initial proximity}

Now that we have constructed approximations of our distributions that avoid recollisions, we can at last compare both the BBGKY and limiting distributions thanks to the coupled pseudo-trajectories. Indeed, we can write the coupling
\begin{align*} 
\widetilde{F}_n^{\varepsilon, [\bV, \delta]} - \widetilde{G}_n^{ [\bV, \delta]} =  & \sum_{\uj_K} \sum_{\uchi_{J_K}} p_\mu^{|\ul^*_{J_K}|}  \int_{T^{[\delta]}_{\underline{j}_K}(t)} \dd \ut_{J_K}  \int_{ \Pi(\uz_n^{[t]}, \uj_K, \uchi_{J_K})^c} \dd \underline{\omega}_{J_K} \dd v_{n+1} \dots \dd v_{N_K} \\ 
& \times \prod_{i=1}^{{J_K}} s_i \langle \omega_{i},  v_{n+i} - v_{m_i}^{[t_i^+]} \rangle_+  \left[ F_{N_K}^\varepsilon \bigl(0, \uz_{N_K}^{[0]} \bigr) - G_{N_K} \bigl(0, \uze_{N_K}^{[0]} \bigr) \right]\ind_{\left\|\uv_{N_K}^{[0]}\right\|^2 \leqslant \bV^2}, \nonumber
\end{align*}
with the same collision parameters for both pseudo-trajectories, the hard-sphere one and the limit one. 
Since by construction $\ux_n^{[t]} =  \uy_n^{[t]}$, and by~Lemma~\ref{lemm:geom} for $\uz_n^{[t]} \in \mE^t_n(\varepsilon_d)$, the $n$ first particles also have identical velocities on~$[0,t]$, then for all times~$\tau \in [0,t]$, we have~$\uz_n^{[\tau]} =  \uze_n^{[\tau]}$ also in positions. 
Henceforth, since by the work done in Section~\ref{sec:disc_labelled} all the added particles are particles at equilibrium tagged~0, one has (as $\ul_n \subset \lbr 1 , n \rbr$),
\begin{align}
G_{N_K} \bigl(0, \uze_{N_K}^{[0]} \bigr)& = M_\beta^{\otimes N_K}(\uu_{N_K}^{[0]}) \varphi_0^{\otimes \ul_n}(\uze_{\ul_n}^{[0]}) \nonumber \\
& = M_\beta^{\otimes N_K}(\uv_{N_K}^{[0]}) \varphi_0^{\otimes \ul_n}(\uz_{\ul_n}^{[0]}) , \label{lgn:traj_comp}
\end{align}
and so 
\begin{align*}
\left| F_{N_K}^\varepsilon \bigl(0, \uz_{N_K}^{[0]} \bigr) - G_{N_K} \bigl(0, \uze_{N_K}^{[0]} \bigr) \right| 
& = \left| F_{N_K}^\varepsilon \bigl(0, \uz_{N_K}^{[0]} \bigr) - G_{N_K} \bigl(0, \uz_{N_K}^{[0]} \bigr) \right| \\
& \leq (C_0 )^{N_K} e^{- \beta \|\uv_{N_K}^{[0]}\|^2} \cdot \varepsilon
\end{align*}
by the initial proximity result~\eqref{ineq:initial_proximity}.
Hence, the previous estimates on the successive-collision operators eventually yield the following result.
\begin{lemm}[Initial value error] Conditionally to $\uz_n \in \mE^t_n(\varepsilon_d)$, the initial error is bounded by
\begin{align}  \label{ineq:coupling}
\left\| \widetilde{F}_n^{\varepsilon, [\bV, \delta]} - \widetilde{G}_n^{ [\bV, \delta]} \right\|_{\Lp^\infty(\mE^t_n(\varepsilon_d))} \leq C^{nK + A^K} \varepsilon.
\end{align}
\end{lemm}

\paragraph{Coherent choice of truncation parameters~$(K, \bV, \delta)$}

Eventually, the last step is to tune the truncation parameters according to~$\varepsilon$, so as to obtain the convergence we want. Actually, there is room to choose the scaling, since one can set all the errors as powers of~$\varepsilon$, except for the pruning error~\eqref{ineq:est_pruned_out} which is significantly bigger. Explicitly, stacking all the errors (pruning~\eqref{ineq:est_pruned_out}, removing additional tagged particle~\eqref{ineq:discard_tagged}, energy truncation~\eqref{ineq:energy_error}, time separation~\eqref{ineq:time_separation} and removing pathological trajectories~\eqref{ineq:path_traj}) then using the coupling result~\eqref{ineq:coupling}, one has
\begin{align*}
&\| F_n^\varepsilon(t) - G_n(t) \|_{\Lp^\infty(\mE^t_n(\varepsilon_d))} \leq C^{C_0\lambda} n^{cn} e^{-2^{K - K^\alpha}} \\
& \hspace{20mm}+  C^{nK + A^K} \Biggl( p_\mu + e^{-\frac{\beta}{2} \bV^2} + \delta +  \Biggl[ \varepsilon^{\frac{d}{d+2}} + \bV^d \times \varepsilon^{\frac{d-1}{3(d+1)}} +  \frac{\bV^\frac{d+1}{2} \varepsilon^{\frac{d(d-1)}{2(d+1)}}}{\delta^{\frac{d-1}{2}}}  \Biggr] + \varepsilon \Biggr)\nonumber.
\end{align*}
Hence, choosing the scaling
\begin{equation*}
 \delta = \varepsilon^\frac{d-1}{d+1},\ \  \und \bV^2 = \frac{2}{\beta}|\log \varepsilon|,
\end{equation*}
one gets
\begin{align}
e^{-\frac{\beta}{2} \bV^2} + \delta +  \bV^d \times \varepsilon^{\frac{d-1}{4(d+1)}} +  \frac{\bV^\frac{d+1}{2} \varepsilon^{\frac{d(d-1)}{2(d+1)}}}{\delta^\frac{d-1}{2}} 
& \leq \varepsilon + \varepsilon^\frac{d-1}{d+1} + \left|\frac{2}{\beta}\log\varepsilon\right|^{\frac{d}{2}} \varepsilon^\frac{d-1}{3(d+1)} + \left|\frac{2}{\beta}\log\varepsilon\right|^{\frac{d+1}{4}} \varepsilon^\frac{d-1}{2(d+1)} \nonumber \\
& \leq \varepsilon^\frac{d-1}{4(d+1)} \label{lgn:vareps_concl}
\end{align}
for $\varepsilon$ small enough. This way, if we pick
\begin{equation*}
K = \left\lfloor \frac{1}{\log A} \log \left( \frac{ (d-1) |\log \varepsilon | }{8(d+1) \log C} \right) \right\rfloor,
\end{equation*}
then  $ C^{A^K} \leq \varepsilon^{- \frac{d-1}{8(d+1)}}$, and we can deal with the term~$C^{nK}$ like in the proof of Proposition~\ref{prop:estrem}, yielding the same factor $n^{cn}$. Hence, for $K$ large enough, denoting $c_\beta \doteq (d-1)/8(d+1)\log C$, we~have
\begin{equation*}
e^{-2^{K - K^\alpha}} \leq e^{-2^{(1-\alpha)K}} \leq \exp\left(- c_\beta |\log \varepsilon|^{(1-\alpha)\frac{\log 2}{\log A}} \right).
\end{equation*}
As $A > 2$ is arbitrary, let us choose $A$ such that $(1-\alpha) \frac{\log 2}{\log A} \geq 1 - 2\alpha$. This gives us the final condition on
\begin{equation*}
\lambda \leq \frac{c_\beta}{2 C_0 \log C} |\log \varepsilon|^{1 - 2\alpha}.
\end{equation*}
Note that we take this scaling so as to have the biggest~$\lambda$ possible in the hypothesis, pushing the associated error at the same level as the truncation error above, which is the limiting one.
It yields
\begin{equation*}
C^{C_0 \lambda} e^{-2^{K - K^\alpha}} \leq \exp\left(- \frac{c_\beta}{2} |\log \varepsilon|^{1-2\alpha} \right).
\end{equation*}
The final verification is the condition $p_\mu \leq 2^{-K} $ of Proposition~\ref{prop:discard}, which is satisfied with room to spare considering the choices above.
Piling all of this, the pruning error is bigger than all the other ones, and we end up with the very last inequality for $\varepsilon$ small enough
\begin{equation*}
 \| F_n^\varepsilon(t) - G_n(t) \|_{\Lp^\infty(\mE^t_n(\varepsilon_d))} \leq n^{cn} \exp\left(- \frac{c_\beta}{2} |\log \varepsilon|^{1-2\alpha} \right).
\end{equation*}
Eventually, the indicator~$\ind_{\mE_n^t(\varepsilon_d)}$ of the set of past-excluding configurations~\eqref{def:past-exclu} pointwise converges, as $\varepsilon$ goes to 0, to the indicator of the following set of full measure
\begin{equation*}
\Bigl\{ \uz_n \in \mD^n\ , \ \forall 1\leq i < j \leq n, \ x_i \neq x_j \und v_i - v_j \notin \Vect(x_i - x_j) \Bigr\},
\end{equation*}
which concludes the proof of Theorem~\ref{theo1}.
\CQFD

\section*{\Huge Appendix}
\setcounter{section}{0}
\renewcommand{\thesection}{\Alph{section}}
\renewcommand{\theHsection}{appendixsection.\Alph{section}}

\section{Cumulants of the exclusion}  \label{app:cumulants_exclusion}

This section is devoted to the study of the cumulants of the exclusion indicator associated with the hard sphere domain~\eqref{def:dom} (see the companion paper~\cite{fou26cumulants} for a generalized definition of the cumulants). 

Denoting $\mG_S$ the set of graphs on a set $S$ and $\mG_n \doteq \mG_{\lbr 1, n\rbr}$ the set of graphs on $\{1, \dots, n\}$, one may write
\begin{align*}
\ind_{\mX_n^\varepsilon}(x_1, \dots, x_n) & = \prod_{1\leq i < j \leq n} (1 - \ind_{x_i \sim x_j}) \\
&= \sum_{G\in \mG_n} \prod_{\{i,j\}\in E_G} (-\ind_{x_i \sim x_j}).
\end{align*}
One can even further decompose the graphs into their connected components. We denote $\mC_S \subset \mG_S$ the set of \emph{connected} graphs on $S$ and $\mP_n$ the set of partitions of $\{1, \dots, n\}$, so that we can write
\begin{align*}
\ind_{\mX_n^\varepsilon}(x_1, \dots, x_n) = \sum_{\sigma \in \mP_n} \prod_{k=1}^{|\sigma|} \left( \sum_{G_k \in \mC_{\sigma_k}} \prod_{\{i,j\}\in E_{G_k}} (-\ind_{x_i \sim x_j}) \right).
\end{align*}
We thus define the cumulants of the exclusion as
\begin{equation}
\phi_k(x_1, \dots, x_k) \doteq \sum_{G \in \mC_k} \prod_{\{i,j\}\in E_G} (-\ind_{x_i \sim x_j}). \label{def:cumulants_of_the_exclusion}
\end{equation}
The very specific structure of these cumulants yields a strong bound on them, called the tree inequality, exposed in the following Proposition. This bound on the integral of the cumulants allows a strong control of the particle correlations, used in Appendix~\ref{app:partition_function} to gain estimates on the partition function.
\begin{prop}[Tree inequality]\ \\ \label{prop:TreeInequality}
\begin{enumerate}[label=(\roman*)] 
\item The modulus of the cumulants may be controlled restricting the sum defining them to the trees~$\mT_k \subset \mC_k$ (i.e. to the minimally connected graphs) as such:
\begin{equation*}
|\phi_k(x_1, \dots, x_k)| \leq \sum_{T \in \mT_k} \prod_{\{i,j\}\in E_T} \ind_{x_i \sim x_j}.
\end{equation*}
\item As a consequence, we have the following control over their integral
\begin{equation*}
\int |\phi_k(\ux_k)| \dd \ux_k \leq k^{k-2} (|\mB_d|\varepsilon^d)^{k-1}.
\end{equation*}
\end{enumerate}
 \end{prop}
\begin{prof} The proof of this proposition relies on a partition scheme due to Penrose~\cite{1967penrose}, and may also be found in~\cite{2023grandev}.
\begin{enumerate}
\item The key argument is to find a map $\pi: \mC_k \rightarrow \mT_k $ such that for any tree $T\in \mT_k$, there is a connected graph $R(T)\in \mC_k$ satisfying
\begin{equation*}
\pi^{-1}(T) = \{ G \in \mC_k \ , \ E_T \subset E_G \subset E_{R(T)} \}.
\end{equation*}
This means that we can partition $\mC_k$ into subsets corresponding each to a single tree, and containing all the graphs that are both compatible with this tree and smaller than an upper graph~$R(T)$. One can find the construction of this partition scheme for example in~\cite{fouthese}.

Now, we decompose the sum defining the~$k$-th cumulant according to this mapping, and using its structure we get
\begin{align*}
\sum_{G\in \mC_k} \prod_{\{i,j\} \in E_G} (-\ind_{i\sim j} )& = \sum_{T\in \mT_k} \sum_{G\in \pi^{-1}(T)} \prod_{\{i,j\} \in E_G} (-\ind_{i\sim j} )\\
& = \sum_{T\in \mT_k} \left(  \prod_{\{i,j\} \in E_T} (-\ind_{i\sim j} ) \right) \left( \sum_{E' \subset E_{R(T)}\setminus E_T} \prod_{\{i,j\} \in E'} (-\ind_{i\sim j} ) \right) \\
& = \sum_{T\in \mT_k} \left(  \prod_{\{i,j\} \in E_T} (-\ind_{i\sim j} ) \right) \left( \prod_{\{i,j\} \in E_{R(T)}\setminus E_T} (1 -\ind_{i\sim j} ) \right) 
\end{align*}
reversing the usual computation so as to harness the cancellations, which yields the result since $1 -\ind_{i\sim j} \in \{0,1\}$.

\item
By the first part of the proposition, the problem is reduced to the integration over the vertices of a tree, for which we know a very simple algorithm. We start integrating over a leaf $i_\ell$ of $T$, which appears in only one condition $\ind_{i_\ell \sim j}$, and we then iterate for a new leaf until the tree is reduced to a single root, so that
\begin{align*}
\int |\phi_k| & \leq \int \sum_{T \in \mT_k} \prod_{\{i,j\} \in E_T} \ind_{i \sim j} \dd \ux_k \\
& \leq  \sum_{T \in \mT_k} |\mB_d|\varepsilon^d \int  \prod_{\{i,j\} \in E_{T\setminus\{ i_\ell \}}} \ind_{i \sim j} \dd x_1 \dots \check{\dd x_{i_\ell}} \dots \dd x_k \\
& \leq \sum_{T \in \mT_k} (|\mB_d|\varepsilon^d)^{k-1} \\
& \ \  = k^{k-2} (|\mB_d|\varepsilon^d)^{k-1},
\end{align*}
concluding with Cayley's formula, that gives the number of trees on $\lbr 1, k\rbr$. \CQFD
\end{enumerate}
\end{prof}

\section{Partition function estimates} \label{app:partition_function}

This section is dedicated to the proof of Proposition~\ref{prop:ZestimeeGC}, which yields the following estimate on the partition function, for $\mu$ large enough in our mixed low density scaling~\eqref{eq:scaling}:
\begin{equation*}
\frac{1}{\mZ_{\mu}}\sum_{q,r\geqslant 0} \frac{ (\lambda C_0)^q \mu^r}{q!r!} \int \ind_{\mX^\varepsilon_{q+r}} \leq C_d^{C_0 \lambda}. 
\end{equation*}
First of all, let us use the symmetry of the particles (among each labelled group) to observe that the partition function can be rewritten as 
\begin{align}
\mZ_{\mu} & = \sum_{p\geqslant 0} \sum_{\ul_p \in \Lambda_p } \frac{\lambda^{|\ul_{p}|} \mu^{p - |\ul_{p}|} }{p!} \int (M_\beta\varphi_0)^{\otimes \ul_{p} }(\uz_{\ul_{p}}) \ind_{\mX^\varepsilon_{p}}(\ux_{p}) \dd \ux_{p} \dd \uv_{\ul_p} \nonumber \\
&= \sum_{p\geqslant 0} \sum_{k=0}^p \frac{\lambda^{k} \mu^{p - k} }{k!(p-k)!} \int (M_\beta\varphi_0)^{\otimes k }(\uz_k) \ind_{\mX^\varepsilon_{p}}(\ux_{p}) \dd \ux_{p} \dd \uv_{k} \label{Zrewrite2}\\
&= \sum_{p,q\geqslant 0} \frac{\lambda^{q} \mu^{p} }{q!p!} \int (M_\beta\varphi_0)^{\otimes q }(\ux_{q}) \ind_{\mX^\varepsilon_{p+q}}(\ux_{p+q}) \dd \ux_{p+q} \dd \uv_{q}. \label{Zrewrite}
\end{align}
We now start by proving the computational lemma below, giving an explicit and exact formulation of the partition function using the cumulants of the exclusion, defined in the Appendix~\ref{app:cumulants_exclusion} above.
\begin{lemm} \label{lemme:Zegal} For all $\lambda, \mu > 0$, with $(\phi_{k})$ the cumulants associated to the exclusion~$(\ind_{\mX^\varepsilon_k})$, one has 
\begin{equation}
\mathcal{Z}_{\mu}  = \exp\left(  \sum_{ (k, l)  \neq (0,0)} \frac{\mu^{k}\lambda^{l}}{k!\ l! } \int (M_\beta\varphi_0)^{\otimes l } \phi_{k+l} \right). \label{eq:Z2esp}
\end{equation}
\end{lemm}

\begin{prof}
We start from the formulation~\eqref{Zrewrite} of the partition function, and we will write the following cumulant expansion, denoting $\mP^{s}_{p+q}$ the partitions of~$\lbr 1, p+q \rbr$ into $s$~subsets,
\begin{align*}
\mathcal{Z}_\mu  & = 1 + \sum_{p,q \neq 0,0} \frac{\mu^p \lambda^q}{p!q!} \sum_{s=1}^{p+q} \sum_{\sigma \in \mathcal{P}_{p+q}^s} \int \prod_{k=1}^{s}  \phi_{|\sigma_k|}(\ux_{\sigma_k}) (M_\beta\varphi_0)^{\otimes q }(\ux_{q}) \dd \ux_{p+q} \\
& = 1 + \sum_{p,q \neq 0,0} \frac{\mu^p \lambda^q}{p!q!} \sum_{s=1}^{p+q} \sum_{ \substack{ (k_i, l_i)  \neq (0,0)  \\  \sum k_j = p, \sum l_j = q }} \frac{1}{s!} \frac{p!}{k_1!\cdots k_s!}\frac{q!}{l_1!\cdots l_s!} \int \prod_{i=1}^{s} \phi_{k_i+l_i} \ (M_\beta\varphi_0)^{\otimes l_i }. \nonumber
\end{align*} 
To perform this last equality, denoting $P = \lbr 1, p \rbr \und Q = \lbr p+1, p+q \rbr$, we consider the following surjection that counts the number of elements from each partition subset in both~$P$ and~$Q$, defined up to an arbitrary choice of a partition order:
\begin{equation*}
 \Phi : \begin{array}{ccc}
\mP^{s}_{p+q} & \longrightarrow  & \Bigl\{ (\uk_s, \ul_s) \ \Bigl| \  (k_i, l_i)  \neq (0,0), \textstyle  \sum k_j = p, \sum l_j = q \Bigr\} \\
\sigma & \mapsto & (|\sigma_1 \cap P|, \dots, |\sigma_s \cap P|, |\sigma_1 \cap Q|, \dots, |\sigma_s \cap Q|).
\end{array}
\end{equation*}
Its defect of injectivity is indeed given by
\begin{equation*}
\sharp \left\{ \sigma \in \mP^{s}_{p+q} \ , \ \Phi(\sigma) = (\uk_s, \ul_s) \right\} = \frac{1}{s!} \cdot \frac{p!}{k_1!\cdots k_s!}\cdot \frac{q!}{l_1!\cdots l_s!},
\end{equation*}
providing the combinatorial factors. Now, permuting the sums so as to get rid of the condition on the value of the index sums, we get
\begin{align}
\mathcal{Z}_\mu
& = 1 + \sum_{s\geqslant 1}^{} \frac{1}{s!} \prod_{i=1}^s \left( \sum_{  (k_i, l_i ) \neq (0,0)}  \frac{\mu^{k_i}}{k_i!}\frac{\lambda^{l_i}}{l_i!} \int \phi_{k_i+l_i}(\ux_{k_i+l_i}) (M_\beta\varphi_0)^{\otimes l_i }(\ux_{l_i}) \dd \ux_{k_i + l_i} \right) \label{ligne:permutesommes}  \\
& = \exp\left(  \sum_{ (k_i, l_i)  \neq (0,0)} \frac{\mu^{k_i}\lambda^{l_i}}{k_i!l_i! } \int \phi_{k_i+l_i} \cdot  (M_\beta\varphi_0)^{\otimes l_i } \right), \nonumber
\end{align}
concluding the proof. \CQFD
\end{prof}

\begin{propf}
Now, Lemma~\ref{lemme:Zegal} applied to both the numerator and the denominator---that have the same structure---gives 
\begin{align*}
& \frac{1}{\mZ_{\mu}}\sum_{p,q\geqslant 0} \frac{\mu^p (\lambda C_0)^q}{p!q!} \int \ind_{\mX^\varepsilon_{p+q}}  = \exp\left(  \sum_{ (k, l)  \neq 0,0} \frac{\mu^{k}\lambda^{l}}{k!\ l! }  \int \phi_{k+l} \cdot  \left[ C_0^l - (M_\beta \varphi_0)^{\otimes l }\right] \right).
\end{align*}
The terms in the sum are vanishing for $l = 0$, so that 
\begin{align*}
\left| \sum_{ (k, l)  \neq 0,0} \frac{\mu^{k}\lambda^{l}}{k!\ l! }  \int \phi_{k+l} \cdot  \left[C_0^l - (M_\beta \varphi_0)^{\otimes l }\right] \right| & \leq \sum_{ \substack{ k\geqslant  0\\ l \geqslant 1}} \frac{\mu^{k}\lambda^{l}}{k!\ l! }  \int |\phi_{k+l}| \cdot  2  C_0^l \\
& \leq 2 \sum_{r \geqslant 1} \sum_{l=1}^r \frac{ \mu^{r-l} (\lambda C_0)^l }{(r-l)! \ l!} \int |\phi_r|\\
& \leq 2 \sum_{r \geqslant 1} \sum_{l=1}^r \frac{ \mu^{r-l} (\lambda C_0)^l }{(r-l)! \ l!} r^{r-2} (|\mB_d| \varepsilon^d)^{r-1},
\end{align*}
by the point \emph{(ii)} of Proposition~\ref{prop:TreeInequality}, stemming from the tree inequality.
Using Stirling's approximation, we have for any $l \leq r$
\[ \frac{r^{r-2}}{(r-l)!\ l!} \leq \frac{e^r}{r^2} \cdot \frac{r!}{(r-l)!\ l!} \leq \frac{(2e)^r}{r^2},\] 
so that
\begin{align*}
\left| \sum_{ (k, l)  \neq 0,0} \frac{\mu^{k}\lambda^{l}}{k!\ l! }  \int \phi_{k+l} \cdot  \left[ C_0^l - (M_\beta\varphi_0)^{\otimes l } \right] \right| & \leq 4e C_0 \lambda \ \sum_{r \geqslant 1} \sum_{l=1}^r \frac{ \mu^{r-l} (\lambda C_0)^{l-1} }{r^2} (2e|\mB_d| \varepsilon^d)^{r-1}.
\end{align*}
Eventually, thanks to the scaling $\mu \varepsilon^{d-1} = 1$ and $1 \ll \lambda \ll \mu$, we have
\begin{align*}
\sum_{r \geqslant 1} \sum_{l=1}^r \frac{ \mu^{r-l} (\lambda C_0)^{l-1} }{r^2} (2e|\mB_d| \varepsilon^d)^{r-1} & \leq \sum_{r \geqslant 1} \sum_{l=1}^r \frac{ \mu^{r-1} }{r^2} (2e|\mB_d| \varepsilon^d)^{r-1} \\
& \leq \sum_{r \geqslant 1} (2e|\mB_d| \varepsilon)^{r-1} \\
& \leq 2,
\end{align*}
concluding the proof. \CQFD
\end{propf} 

%%%--%%%--%%%--%%%--%%%-- Bibliographie --%%%--%%%--%%%--%%%

\newpage

\bibliographystyle{abbrv}  % abbrv  acm   alpha (noms complets)    apalike   ieeetr ...
%\bibliography{../../arefs}

\bibliography{arefs}

\end{document}